\documentclass[11
pt]{article}

\newcounter{alfej}[section]

\newcommand{\bekezdes}[2]{\vspace*{-2mm}\paragraph{#1}\label{#2}}

\usepackage{amssymb}

\begin{document}

\let\d\partial
\def\EE{\mathcal E}
\def\C{\mathbb C}
\def\Q{\mathbb Q}
\def\R{\mathbb R}
\def\Z{\mathbb Z}
\def\Res{\operatorname{Res}}
\def\im{\operatorname{im}}
\def\Ker{\operatorname{Ker}}
\let\ker\Ker
\def\Tot{\operatorname{Tot}}
\def\Gr{\operatorname{Gr}}
\def\Hom{\operatorname{Hom}}
\def\Log{\operatorname{Hom}}
\def\codim{\operatorname{codim}}
\let\emptyset\varnothing
\def\resp{\operatorname{respectively}}
\def\resp{\operatorname{respectively}}
\def\mod{\operatorname{mod}}
\def\modulo{\operatorname{modulo}}
\def\b{\bullet}
\def\c{\mathcal}
\def\p{\pitchfork}
\def\Y{Y_{\alpha_1,\ldots,\alpha_p}}
\def\u{\underline}
\def\y{\tilde{Y}^p}
\def\s{\star}
\def\noi{\noindent}
\newcommand{\be}{\bekezdes} 
\def\esq{$\hfill\square$}

\begin{center}
{\Large \bf Hodge-DeRham theory with degenerating coefficients}\\
 \ \\ {\Large Fouad Elzein}
\end{center}

\begin{abstract}
 Let ${\cal L}$ be a local system  on the complement $X^{\star}$ of
 a normal crossing divisor (NCD) $ Y$ in a smooth analytic variety $X$
and let $ j: X^{\star} = X - Y \to X $ denotes the open embedding. The purpose of this 
paper is to describe a weight filtration $W$ on the direct image ${\bf j}_{\star}{\cal L} 
$  and in case 
a morphism $f: X \to D$ to a complex disc is given  with $Y = f^{-1}(0)$, the weight 
filtration on 
the complex of nearby cocycles $\Psi_f ({\cal L})$ on $Y$. A comparison theorem shows 
that the filtration coincides with the weight defined by the 
 logarithm of the monodromy and provides the link with various results on the subject.

\vspace{3mm}

(secondary).

\end{abstract}

\section{Introduction}
We consider a  local system ${\cal L}$ on the complement $X^{\star}$ of
 a normal crossing divisor (NCD) $ Y$ in a smooth analytic variety $X$
and let $ j: X^{\star} = X - Y \to X $ denotes the open embedding. The purpose of this 
paper is to describe a weight filtration $W$ on the direct image ${\bf j}_{\star}{\cal L} 
$ ( denoted  also $ {\bf R}j_{\star}{\cal L} $)
in the  category of perverse sheaves  and in case 
a morphism $f: X \to D$ to a complex disc is given with $Y = f^{-1}(0)$, the weight 
filtration on 
the complex of nearby cocycles $\Psi_f ({\cal L})$ on $Y$.
This subject started with Deligne's paper [8] when
 ${\cal L} = \ C$ ( Steenbrink treated the complex  of nearby cycles in [32]), then 
developed 
extensively again after the discovery of the theory of intersection
cohomology [17] perverse sheaves [2] and 
the purity theorem [2], [5],
[24]. The theory of differential modules introduced by Kashiwara proved to be fundamental 
in the understanding 
of the problem as it appeared in the work of Malgrange [28], Kashiwara [24]-[26] and 
later M. Saito [29],[30]. 
It is interesting
to treat the problem by the original logarithmic methods as proposed in the note [14], in 
order to obtain 
topological interpretation of these results in the direction of [27](for another 
direction [15]).
Let me  explain now the problems encountered in the construction of a mixed Hodge 
structure
(MHS) for $X$ proper  on the cohomology $ H^*( X-Y, {\cal L})$.  

\smallskip

\noindent{\it  Hypothesis 1}. Let  ${\cal L}$ be a local system defined
over $\Q$, on the complement of the normal crossing divisor (NCD) $ Y$
in a smooth analytic variety $X$, $({\cal L}_X,
\nabla) $ the canonical extension  of ${\cal L}^{\C} = {\cal L} \otimes \C $ [6] with a 
meromorphic connection $\nabla$ having a  regular singularity
along  $ Y$ in $X$ and the associated DeRham 
logarithmic complex $ {\Omega}^*_X (Log Y) \otimes {\cal L}_X $ defined by $\nabla$ (in 
the text we write  $ L $ and ${\cal L}  $ for the rational as well complex vector spaces 
).

\smallskip

In order to construct a weight filtration $ W $ by subcomplexes of 
$ {\Omega}^*_X (Log Y) \otimes {\cal L}_X $ we need a precise
description of the 
correspondence with the local system $ {\cal L}$. For all subset $M$ of
$I$, let 
$Y_M = {\cap_{i \in M}} Y_i$ and $Y_M^* = Y_M - {\cup_{i \not\in M}}
Y_i$, $j_M 
\colon Y^*_M \rightarrow Y_M$  the locally closed
embeddings, then $Y_M - Y_M^*$ is a NCD in $Y_M$ and the open subsets $Y_M^*$  of $Y_M$
form with $X^*$ a natural stratification of $X$ ( we suppose the NCD $Y = \cup_{i \in I} 
Y_i $  equal to the union of irreducible  and smooth components $Y_i$ 
for $i$ in  $I$). All extensions of  $\cal L$ that we will introduce will be 
constructible with respect to this stratification and even perverse. In fact we will need 
a combinatorial model of $ {\Omega}^*_X (Log Y) 
\otimes {\cal L}_X $ for the description of the weight.\par

If we consider a point $y \in Y_M^*$ and a variation of Hodge structures on 
${\cal L} $ of weight $m$, locally defined by a nilpotent orbit $L$ and a set of 
nilpotent endomorphisms $N_i, i \in M $, the nilpotent orbit theorem [4], [24]
states that it degenerates along $Y_M^*$ into a variation of $ MHS $ with weight 
filtration $W^M = W(\Sigma_{i\in M} N_i)$ shifted by $m$, however this result doesn't 
lead directly to a structure of mixed Hodge complex since what happens at the 
intersection of $Y_M$ and $Y_K$ for two subsets $M$ and $ K$ of $I$ couldn't be explained 
until the discovery of perverse sheaves. The $MHS$ we are 
looking for cannot be obtained from  Hodge complexes defined by smooth and  proper 
varieties, so it was only after the purity theorem [2], [5], [24]
and the work on perverse sheaves [2],
that  the weight filtration $\cal {W} (\cal {N})$ on $\Psi_f ({\cal L})$  has been  
defined by the logarithm
of the monodromy  $\cal {N}$ in the abelian category of perverse sheaves ( object not 
trivial to compute). In characteristic zero this has been
successfully related to 
the theory of differential modules [29]. In this work we study the construction of the 
weight filtration given in the note [14] with new
general proofs of the purity and decomposition  statements  there. The
key result that enables us to simplify most of 
the proofs is the following decomposition (see \S 2 ):

$$Gr^{W^M}_r  L  = \oplus_{\Sigma_{i\in M} m_i = r} Gr^{W^n}_{m_n} \cdots
Gr^{W^1}_{m_1} L $$

\noindent and in fact the intersection complexes defined by the local systems on $Y^*_M$ 
with fiber:

$$\oplus_{\Sigma_{i\in M} m_i = r - \vert M \vert, m_i \geq 0
} Gr^{W^n}_{m_n} 
\cdots  Gr^{W^1}_{m_1} (L / (N_1 L + \cdots + N_n L)) $$

\noindent will be the Hodge components of the decomposition of the graded part of the 
weight filtrations (of the  
primitive parts in the case of nearby cocycles). This main result in the open 
case, its extension to the case of nearby cocycles   and its comparison with 
Kashiwara and Saito's results form the
contents of the article in the second and third sections. We hope that this 
comparison will be helpful to the reader who doesn't want to go  immediately through
the whole subject of differential modules.

One important improvement in the new purity theory with respect to the theory of mixed 
Hodge complexes in [8], [9] is 
the fact that the objects  can be defined locally: these objects are the  
intermediate extensions of $VHS$ on $Y^*_M $ for some $M\subset I$ since such extension 
is defined on $Y_M $ even if it is not proper which makes such 
objects very pleasant to use especially in the $\Psi_f$ case. We adopt here the 
convention to use the terminology of perverse sheaves up to a shift in degrees, although 
it is very important to give explicitly this shift when needed in a 
proof. Let us review the contents: the definition of the weight is in $(\S 1, 
II,2)$ formula $(14)$; the main results are in $\S 2$, namely the key lemma for 
the decomposition at the level of the weights $(I,1)$, the local purity 
$(I.3)$, the local decomposition $(II)$ and the global decomposition $(III)$.
 In $\S 3$, the weight of the nilpotent action on $\Psi_f$ is in $(I)$, the  
local decomposition in $(I,1)$, the global decomposition in $(I,2)$, the comparison 
 in $(II)$ and in an  example we apply this theory to remove the base change in 
Steenbrink's work. In the last part we just state the results for good $VMHS$
as there is no new difficulties in the proofs. Finally we suggest strongly to the reader 
to follow the proofs on an example, sometimes on the surface case as in $(II.1.5)$, or 
for $X$ a line and $Y = 0$ a point, then the fiber at $0$ of ${\bf j}_{\star}{\cal L} $ 
is a complex $ L \stackrel{N}{\rightarrow} L$ where $(L, N)$ is the nilpotent orbit   of 
weight $m$ defined at $y$ and the weight on the complex  is ${\cal W}[m]$ defined by
${\cal W}[m]_{r+m} = {\cal W}_r = (W_{r+1}L \stackrel{N}{\rightarrow} W_{r-1}L)$. This 
example will be again useful for $\Psi_f {\cal L}$ in $\S 3$.

\subsection{Preliminaries on perverse extensions and nilpotent orbits}

 In the neighbourhood of a point $y$ in 
$Y$, we can suppose $X \simeq D^{n+k}$ and $X^* \simeq {(D^*)}^n \times D^k$ 
where $D$ is a complex disc, denoted with a star when the origin is deleted. The 
fundamental group $\Pi_1 (X^*)$ is a free abelian group generated by $n$ elements 
representing  classes of  closed paths around the origin, one for  each 
$D^* $ in the various axis with one dimensional  coordinate $z_i$  ( the hypersurface 
$Y_i$ is defined by the equation $z_i=0$ ) . Then the local
system $\cal L$ corresponds to a representation of $ 
{\Pi}_1(X^*)$ in a vector space  $L$, i.e the action of commuting automorphisms 
$T_i$ for $i \in [1,n ]$ indexed by the local components $Y_i$ of $Y$ and called 
monodromy action around $Y_i$. The automorphisms $T_i$ decomposes as a product 
of commuting automorphisms,  semi-simple and unipotent $T_i = T_i^s T_i^u$. When 
$L$ is a $\C $ - vector space, $T_i^s$ can be represented by the diagonal matrix of its 
eigenvalues. If we consider sequences of eigenvalues $a_i$ for 
each $T_i$ we have the spectral  decomposition of $L$ 
$$L = {\oplus }_{a .} L^{a .} \qquad , \qquad L^{a .}  = \cap_{i \in [1,n ]} 
({\cup}_{j>0} \; \hbox{ker} \; {(T_i - {a}_i I)}^j) \leqno{(1)}$$

\noindent  where the direct sum is over all families $ (a .) \in \C^n .$  
The logarithm of $T_i$ is  defined as the  sum 
$$Log T_i = - 2i \pi ( D_i + N_i ) =  Log T_i^s + Log T_i^u
\leqno{(2)}$$ where  $ - 2i \pi D_i = Log T_i^s$ is the diagonal matrix
formed by $Log \,{a}_i$  for all
eigenvalues ${a}_i$ of $T_i^s$ and for a fixed determination of Log on
$ {\C}^* $, while $ - 2i \pi N_i = Log T_i^u$ is defined by the
polynomial
function   
$ - 2i \pi  N_i = {\Sigma}_{k \geq 1} (1/ k) {(I-T^u_i)}^k $ in the
nilpotent 
$(I - T^u_i)$ so that the sum is finite. \par

\noindent  Now we  describe various local extensions of $\cal L$.
\be {The (higher) direct image ${\bf j}_{\star}{\cal L} $,
local and global description}{1.I.1} \
  The
complex $\hbox {\bf j}_* \cal L$ is perverse and its fiber at the origin in  
$D^{n+k}$ is quasi-isomorphic to a Koszul complex as follows.
We associate to $(L, (D_i + N_i), i \in [1,n ])$   a strict simplicial
vector space such that for all sequence $(i.) = (i_1 < \cdots < i_p)$ 
  
\centerline {$L(i.) = L \qquad , \qquad D_{i_j} + N_{i_j} \colon L(i. - i_j) 
\rightarrow L(i.)$.}

The associated simple complex is the  Koszul complex or  the exterior 
algebra defined  by $(L, D_i + N_i)$ denoted by $\Omega (L, D. + N.) \colon = 
s(L(J),D.+N.)_{J \subset [1,n ]} $
 where  $J$ is identified with the strictly increasing 
sequence of its elements and where $ L (J) = L $. 
It can be checked that its  cohomology is the same as  $\Omega (L, Id - T.)$, 
the Koszul complex defined by $(L, Id - T_i),  i \in [1,n ]$. 

When we fix a family $\alpha_j \in [0,1 [$ for $ j \in [1,n]$ such that $ e(\alpha_j) = 
e^{- 2i \pi {\alpha}_j} = a_j $ is an eigenvalue for $D_j$, we have
$$\Omega (L, D. + N.)\colon = s(L(J),D.+N.)_{J \subset [1,n ]} = {\oplus}_{\alpha .} 
\Omega (L^{e(\alpha .)},  \alpha . Id + N.).\leqno{(3)}$$ 
In particular each sub-complex is  acyclic when $\alpha _j$ is not zero since 
then $  \alpha _j Id + N_j$ is an isomorphism.
This local setting compares to the global case via Grothendieck and
Deligne DeRham cohomology results. Let $y \in Y^{\star}_M$, then
$$(\hbox {\bf j}_* {\cal L})_y \simeq ({\Omega}^*_X (Log Y) \otimes {\cal L}_X)_y \simeq 
\Omega (L, D_j + N_j, j \in M )\leqno{(4)}$$
 When we start with a $\Q$ local system, ${\cal L}_X$ is equal to the canonical extension 
with residue 
 $\alpha_j \in [0,1[ \cap \Q $ such that $ a_j = e^{- 2i \pi {\alpha}_j}$.
An element $v$ of $L^{e(\alpha .)}_J = L^{e(\alpha .)},  J \subset M,$
corresponds to the section $\mathaccent"7E{v}\otimes \frac {dz_J}{z_J}$ of 
$({\cal L}_X \otimes \frac {dz_J}{z_J})_y$,
where $ \frac {dz_J}{z_J} = \Lambda_{j \in J} \frac {dz_j}{z_j}$   by the formula
$$ \mathaccent"7E{\hbox{v}} (z)  = (exp (\Sigma_{j \in J}  log
z_{j}({\alpha}_{j} + N_{j}))). v =
{\Pi}_{j \in J} {z_{j}}^{{\alpha}_{j}} \, \hbox{exp} ({\Sigma}_{j \in J} log z_{j} 
N_{j}). v  \leqno{(5)}$$ 

\noindent a basis of $L$ is sent on a basis of  $({\cal L}_X)_y$, the  endomorphisms $ 
D_j + N_j $ defines corresponding endomorphisms denoted by the 
same symbols on the image  sections $ \mathaccent"7E{\hbox{v}}$  and we have $$\nabla 
\mathaccent"7E {v} = \Sigma_{j \in J} (\alpha_{j} + N_{j}). 
\mathaccent"7E {v } \otimes \frac {dz_j}{z_j}\leqno{(6)}$$ 
This description of  $(\hbox {\bf j}_* {\cal L })_y$ is the model for the 
description of the next various perverse sheaves.
 \be{The intermediate extension ${\bf j}_{!\star}{\cal L} $}{1.I.2} \
   Let $(D + N)_J = \Pi_{j \in J} (D_j + N_j)$ denotes a composition of endomorphisms of 
$L$, we consider  the strict  simplicial  sub-complex
of the DeRham logarithmic complex (2.4) defined by 
$Im (D + N)_J$ 
in $L(J) = L$. The associated simple complex will be denoted by 
$$IC (L)\colon = s({(D+N)_J L},(D. + N.))_{J \subset [1,n]}
\leqno{(7)}$$
The intermediate extension ${\bf j}_{!*} {\cal L}$ of ${\cal L}$ is defined by 
an explicit formula  in terms of the stratification [24, \S 3]. Locally its  fiber 
at a point $y\in Y^*_M$  is given in terms of the above complex $$ {\bf
j}_{!*}  ({\cal L})_y \simeq IC (L) \simeq s({(D+N)_J L,(D. +
N.))}_{J \subset M}     \leqno{(8)}$$
The corresponding global DeRham description is given as a sub-complex $IC (X, {\cal L})$ 
of 
 ${\Omega}^*_X (Log Y) \otimes {\cal L}_X$. The residue of the
connection 
$\nabla $ along each  $Y_j$ defines an endomorphism 
$ ({\cal D}_j + {\cal N }_j )$ on the restriction 
${\cal L}_{Y_j} $ of ${\cal L}_X$, 
then in terms of a set of $n$ coordinates $y_i, i\in [1,n]$ defining $Y_M$  on an open 
set 
$U_y$ containing $y\in Y^*_M$ where we identify $M$ with $[1,n]= M$ ( $ n = \mid M\mid$) 
and  a section 

\centerline {$ f = \Sigma_{J\subset M,J'\cap M = \emptyset } f_{J,J'} (\frac 
{dy_J}{y_J})\wedge dy_{J'}$ } 
where we define $y_J = \Pi_{j\in J} y_j$ and $dy_J = 
\wedge_{j\in J} d y_j$, we have
$$ f \in IC ( U_y , {\cal L}) \Leftrightarrow \forall J \subset M,\, 
f_{J,J'}/Y_J \in ({\cal D}_J + {\cal N}_{J})({\cal L}_{Y_J})
\leqno{(9)} $$ A global definition of  $IC (X, {\cal L})$, using notations of [26], is 
given 
as follows.  
Consider $M\subset I$ and for all $J \subset M $ the families of sub-bundles 
$({\cal D}_J +{\cal N }_{J})({\cal L}_{Y_J})$ of ${\cal L}_{Y_J}$, then define for each 
$M$
 the sub-module   ${\cal I}(M) $ of ${\cal L}_X$ consisting of sections $v$ satisfying:
 $$ v \in {\cal I}(M) \Leftrightarrow \forall J\subset I, v/Y_J \in ({\cal D}_{M 
\cap J}+{\cal N }_{M \cap J})({\cal L}_{Y_J}) $$ 
This submodule is well defined since we have, for  $J \subset K$, the inclusion 
$({\cal D}_J+{\cal N }_{J})({\cal L}_{Y_J})_{Y_K} \supset ({\cal D}_K + {\cal N 
}_K)({\cal L}_{/Y_K})$. Let $\Omega (M)$ be the ${\cal O}_X $
sub-exterior algebra of ${\Omega}_X^*(Log Y)$ generated by ${\Omega}^1_X(Log Y_i)$ for 
$i$ in 
$M$. Then we can define  
$$IC(X,{\cal L}) = \Sigma_{M \subset
I} \Omega (M) \otimes {\cal I} (M) \subset {\Omega}_X^*(Log Y)\otimes{\cal L}_X.  
\leqno{(10)}$$ 

In terms of the decomposition (3), since   the  endomorphism 
$ ( \alpha_j Id + {\cal N }_j) $  is an isomorphism on ${\cal L}^{e(\alpha .)}_{Y_j} $ 
whenever  $\alpha_j \neq 0 $, we introduce for each set $\alpha .$ 
the subset $I(\alpha .) \subset [1,n ]$ such that $j 
\in I( \alpha .)$ iff ${\alpha}_j = 0$, then for each $J \subset [1,n]\quad , {(D + N)}_J 
L^{e(\alpha .)} = N_{J \cap I( \alpha .)} L^{e(\alpha .)}$ where $N_{J \cap I (\alpha .)} 
= {\Pi_{j \in J \cap I(\alpha .)}} N_j$ ( the identity if $J \cap I (\alpha .)$ is 
empty). 
Then locally the  fiber of ${\bf j}_{!*} {\cal L}$ at a point $y\in Y^*_M$  is 
$$ {\bf j}_{!*}  ({\cal L})_y \cong IC (L) \simeq  \oplus_{\alpha .} s(N_{J \cap 
I (\alpha .)} {L}^{e(\alpha .)})_{J \subset M} \leqno{(11)} $$
       
{\em Remark: The fiber of the complex  $({\Omega}^*_X (Log Y) \otimes {\cal L}^{e (\alpha 
.)}_{X})_y$ is acyclic if there exists an index $j \in M $ such that $\alpha_j \neq 0 $.} 
\be {Hodge filtration, 
Nilpotent orbits and Purity }{1.I.3} \ 
 {\it  Hypothesis 2: Variation of Hodge structures (VHS).} Consider the 
flat bundle $({\cal L}_X, \nabla) $ in the previous hypothesis and suppose now 
that ${\cal L}_{X^*}$ underlies a $VHS$ that is a polarised filtration by 
subbundles $F$ of weight $m$ satisfying Griffith's conditions [19]. 

\smallskip

The {\it nilpotent orbit} theorem [19], [4], [24], [25],  states that $F$ extends 
to a filtration by subbundles $F$ of   ${\cal L}_X$ such that the restrictions 
to open intersections $Y^*_M $ of 
components of $Y$ underly a variation of mixed  Hodge structures $VMHS$ where the 
weight filtration is defined by the nilpotent endomorphism ${\cal N}_M$ defined by the 
connection.  

\smallskip

{\it Local version.}

\smallskip
\noindent Near a point $y \in Y^*_M$ with $\vert M\vert = n$
a neighbourhood of $y$ in the fiber of the normal bundle looks like a disc $D^n$ and the 
above hypothesis reduces to

\smallskip
\noindent {\it Local Hypothesis 2:  Nilpotent orbits [4].} Let 
$$ (L, N_i, F, P, m, i\in M = [1,n] ) \leqno{(12)} $$
be defined by the above hypothesis, that is a $\Q $ vector space $L$ with  endomorphisms 
$N_i$ viewed as defined by the horizontal (zero) sections of the connection on $ 
(D^*)^n$, a Hodge structure $F$ on $L^{\C} = L \otimes_{\Q} \C $ viewed as  the 
fiber of the vector bundle ${\cal L}_X $ at $y$ (here $y = 0 $), a natural 
integer $m$  the weight  and the polarisation $P $. 

\smallskip
The main theorem [4] states that   for all 
$N = \Sigma_{i \in M} {\lambda}_i N_i$ 
with ${\lambda}_i > 0$ in ${\hbox {\bf R}}$ the 
 filtration $W(N)$ ( with center 
$0$ ) is independant of $N$ when  ${\lambda}_i$ vary and $ W(N)[m]$ is the 
weight filtration of a polarised $MHS$ called the  limit $MHS$  of weight $m$ ( $L, 
F,W(N)[m]) $. 
 
\smallskip
\noindent {\em Remark: $ W(N)[m]$ is $ W(N)$ with indices shifted by $m$ to the 
right: $ (W(N)[m])_r\colon= W_{r-m}(N)$, the convention being a shift to left for a 
decreasing 
filtration and to right for an increasing filtration.}

\smallskip
\noindent We say that $ W(N)$ defines a $MHS$ of weight $m$ on $L$. It is very important 
to notice that the same orbit underlies other different orbits
depending on the intersection of components of $Y$ (here the
intersection of the 
axis of $D^n$) where the  point $z$ near $y$ is considered, in particular $F_z \neq F_y 
$. In this case when we restrict the orbit to $J \subset M $, we should write  
   $$ (L, N_i, F (J), P, m, i\in J \subset M ) $$ Finally we
will need the following result [4 p 505]: Let $I, J \subset M$ and let $
N_J  = \Sigma_{i \in J}  N_i$ then $W(N_{I\cup J}) $ is the  weight filtration of $N_J$ 
relative to $W(N_I)$
$$ \forall j,i \geq 0, N_J^i: Gr^{W(N_{I\cup J} )}_{i+j}   Gr^{W(N_I)}_j
{\buildrel \sim \over \rightarrow } Gr^{W(N_{I\cup J} )}_{j-i}
Gr^{W(N_I)}_j.$$

 \subsection{The weight filtration on the logarithmic complex}

Now we want to give the construction of the weight filtration given in
[14] and  based on a general formula of the  intersection  complex given by Kashiwara and 
Kawai[26]. Earlier work in the surface case using ad hoc methods showed that the 
purity and the decomposition theorems could  be obtained out of similar considerations.

To this purpose we introduce a category $S(I)= S$ attached to a set $I$.
We start with a local study, that is to say with the hypothesis of a
polarised 
nilpotent orbit and
we describe the weight filtration $W$ on the DeRham complex $ \Omega (L,N.)$. In fact the 
filtration $W$  is defined on a 
quasi-isomorphic 
complex and may appear unrealistic at first sight, however the features of the 
purity theory will appear relatively quickly. First we ask the reader to take some time 
to get acquainted with the new category $S(I)$ serving as indices for
the new complex. The lowest weight is given by  the intermediate
extension of ${\cal L}$ or $ IC(L)$, then for the higher weights we need to introduce the 
complexes $C^{KM}_r L$ for $ K \subset M \subset I $ which describe the
geometry of the decomposition theorem (\S 2, II) and the purity theory
(\S 2, I.3) where the proof reflects deep relations between
the weifgt filtrations of the various $N_i$.

\be {Complexes with indices in the category $S(I)$.} {1.II.1} \ 
We introduce a category $S(I)= S$ attached to a set $ I$, whose objects consist 
of sequences of  increasing subsets of $I$ of the following form:

$$ (s_.) = (I = s_1 \supsetneqq s_2 \ldots \supsetneqq s_p \neq \emptyset), \; (p > 0)$$

Substracting a subset $s_i$ from a sequence $s_.$ defines a morphism
$\delta_i (s.) : (s. - s_i)\rightarrow s.$ and more generally 
$Hom ( s'., s.)$ is equal
to one element iff $(s.') \leq (s.)$ . We write $s. \in S$ and define its degree $\mid 
s.\mid $ as the number of subsets $s_i$ in $s_.$ or length of the sequence.

\medskip

{\em Correspondence with an open simplex. } If $I = \{1, \ldots, n \}$
is finite, $S(I)$ can be realised as a barycentric subdivision of a simplex of dimension 
$ n - 1 $, a subset $K$ corresponding to the barycenter of
the vertices in $K$ and a sequence of subsets to the simplex defined  by the vertices 
associated to the subsets. Since all sequences contain $I$,
all corresponding simplices must have the barycenter as vertex, that is:
$S(I)$ define a simplicial object computing the $n^{th}$ homology group
with closed support of the open simplex. This remark leads us to the next definition.

\medskip

{\em Simplicial complex defined by complexes with indices in S(I)}. An
algebraic variety over a fixed variety $X$ with indices in $S$ is a
covariant functor $\Pi : X_{s.} \rightarrow X$ for $s.\in S$. An abelian
sheaf (resp. complex of abelian sheaves )${\cal F}$ is a family of
abelian sheaves (resp. complex of abelian
sheaves) ${\cal F}_{s.}$ over  $X_{s.}$ and functorial morphisms ${\cal 
F}_{s.}\rightarrow {\cal F}_{s.'}$ for $(s.') \leq (s.)$. 

The direct image of an abelian sheaf (resp. complex of sheaves ) denoted
$\Pi_*{\cal F}$ or  $s({\cal F}_{s.})_{s.\in S}$ is the simple 
complex(resp.simple complex associated to a  double complex)  on $X$:
$$\Pi_*{\cal F} \colon= \oplus_{s. \in S}(\Pi_*{\cal F}_{s.}) [\mid
s.\mid - \mid I \mid  ], \; d = \Sigma_{i\in [1,\mid s.\mid]} (-1)^i
{\cal F}(\delta_i (s.)).$$ The variety $X$ defines the constant variety
$X_{s.}= X$. The constant sheaf 
$\Z$ lifts to a sheaf on $X_{s.}$ such that the diagonal morphism :
$\Z_X \rightarrow \oplus_{\mid (s.)\mid = \mid I \mid}  \Z_{X_{s.}}$ (
that is : $n \in \Z \rightarrow (\ldots, n_{s.}, \ldots) \in
\oplus_{\mid (s.)\mid = \mid I \mid} 
\Z $ defines a quasi-isomorphism $\Z_X \cong \Pi_*\Pi^* (\Z_X)$. This is
true since $S(I)$ is isomorphic to the category defined by the 
barycentric 
subdivision of an open simplex of dimension  $\mid I \mid - 1$.

\be { Local definition of the weight filtration.}{1.II.2} \
 
\noindent Our hypothesis here consists again of the  nilpotent
orbit $ (L, (N_i)_{i\in M}, F ,m, P ) $   of weight $m$ and polarisation
$P$ and  the corresponding
filtrations $ (W^J)_{J\subset M} $.

\noindent We will use the category $S(M)$ attached to $M $ whose objects
consist of 
sequences of  decreasing subsets of $M$ of the  form $ (s_.) = (M = s_1 
\supsetneqq s_2 \ldots \supsetneqq s_p \neq \emptyset),\;\; p>0.$

\noindent In this construction we will need double complexes, more
precisely complexes of the 
previously defined  exterior complexes. They correspond to objects with
indices in the category $M_.^+\times S(M)$ where the objects of $M_.^+$ are the subsets  
$J\subset M$ including the empty set.
 Geometrically  $M$ corresponds to a normal section to $Y^*_M$ in $X$
and  $J$  to $\wedge_{i\in J}dz_i$ in the exterior DeRham complex written as $s(L_J)_{J 
\subset M}$ on 
the normal section to  $Y^*_M$  and 
the decomposition $M_.^+\simeq 
(M-K)_.^+ \times K_.^+$ corresponds to the isomorphism $\C^M\simeq \C^{(M-K)}\times 
\C^K$.

\smallskip

\noindent {\it Notations.} For each $s.\in S(M)$, $J\subset M$ and integer $r$ we define  
$a_{s_{\lambda}}(J,r) = \mid s_{\lambda}\mid - 2 \mid s_{\lambda}\cap J\mid+ r $, and 
 for all $(J,s.)\in M_.^+\times S(M)$ the vector spaces
$$ W_r(J,s.) \colon = \bigcap_{s_{\lambda}\in
s.}W^{s_{\lambda}}_{a_{s_{\lambda}}(J,r)} 
L,\quad F^r(J,s.)\colon= F^{r-\mid J\mid } L, \; \;
(a_{s_{\lambda}}(J,r) = \mid 
s_{\lambda}\mid - 2 \mid s_{\lambda}\cap J\mid+ r ) $$
where $ W^{s_{\lambda}} $ is centered at $0$, then we define on the DeRham 
complex $\Omega (L,N.)$, the  filtrations by sub-complexes $W(s.)$ (weight) and  $F(s.)$ 
(Hodge) as
$$   W_r(s.)= s(W_r(J,s.))_{J\subset M},\quad F^r(s.)\colon= 
s(F^{r}(J,s.))_{J\subset M}\leqno{(13)}$$
and finally

\noindent {\em Definition: The weight (centered at zero) and Hodge filtrations on  the 
combinatorial  DeRham complex 
 $\Omega^{\star} L = s(\Omega (L,N.))_{s.\in S(M)}$ are defined by ``summing'' over $s.$ 
as:}
  $$F^{r} (\Omega^{\star} L) \colon=  s(F^r(s.))_{s.\in S(M)} \subset 
s(L(s.))_{s.\in S(M)},\;{\cal W}_r (\Omega^{\star} L)\colon=
s(W_r(s.))_{s.\in 
S(M)} \subset s(L(s.))_{s.\in S(M)}  \leqno{(14)}$$

\noindent {\em  Ex in dimension $2$}.

\noindent Let $W^{1,2} = W(N_1 + N_2), W^1 = W(N_1 ) $ and $ W^2 = W(N_2) $, 
the weight ${\cal W}_r $ is a double complex where the first line is the direct sum for $ 
\{1,2 \}\supsetneqq 1 $ and $ \{1,2 \} \supsetneqq 2 $ of:

$W_{r+2}^{1,2}\cap W^1_{r+2}\stackrel{(N_1,\; N_2)} {\longrightarrow}  W_{r}^{1,2}\cap 
W^1_{r}\oplus
W_{r}^{1,2}\cap W^1_{r+2}\stackrel{(- N_2,\;N_1)} \longrightarrow  W_{r-2}^{1,2}\cap 
W^1_{r}$

\noindent and 

$W_{r+2}^{1,2}\cap W^2_{r+2}\stackrel{(N_1,\; N_2)} {\longrightarrow}   W_{r}^{1,2}\cap
W^2_{r+2}\oplus W_{r}^{1,2}\cap W^2_{r} \stackrel{(- N_2,\;N_1)} \longrightarrow  
W_{r-2}^{1,2}\cap
W^2_{r}$  

\noindent The second line for $ \{1,2 \}$ is

$W_{r+2}^{1,2} \stackrel{(N_1,\; N_2)} {\longrightarrow}  W_{r}^{1,2}\oplus 
W_{r}^{1,2}\stackrel{(- N_2,\;N_1)} \longrightarrow
W_{r-2}^{1,2}$.

\noindent which reduces to the formula in [26] for  $ r = m $. 

\be { The Complexes $C_r^{KM}L$ and $ C^K_r L $.}{1.II.3} \

To study the graded part of the weight, we need to introduce the
following subcategories:

\noindent For each subset $ K \subset M$, let 
$ S_K(M) = \{ s. \in S(M): K \in s. \}$ and consider the isomorphism of categories:
$$S(K)\times S(M-K) {\buildrel \sim \over \rightarrow }
 S_K(M), \; (s., s'_.) \rightarrow  ( K \cup s'_., s.)$$
\noindent We consider the vector space with indices $(J,s.) \in M.^+ \times S_KM,$ and 
its associated complex 
 $$C_r(J,s.)L\colon=  \bigcap _{K\neq s_{\lambda}\in s.} 
W^{s_{\lambda}}_{a_{\lambda}(J,r-1)}  Gr^{W^K}_{a_{K(J,r)}} L, \;\;\;C_r (s.)L \colon=  
s(C_r(J,s.)L )_{J \in M^+_.} \leqno{(15)}  $$
we define as well $C_r(J,s.)N_J L $ and $ C_r(J,s.)(L/N_J L)$ by replacing $L$ with 
$N_JL$ and $L/N_JL$, then  
 $ C_r (s.)IC L \colon= s(C_r(J,s.)N_JL )_{J \in M^+_.} $
and $ C_r (s.)QL \colon= s(C_r(J,s.) L/N_JL )_{J \in M^+_.} $

\noindent {\em Definition: For $K \subset  M$ the complex $C_r^{KM} L $ is defined as 
$$ \;C_r^{KM}L\colon=  s(C_r(s.)L)_{s.\in S_K(M)}, \;C_r^{KM}ICL\colon=  
s(C_r(s.)ICL)_{s.\in S_K(M)}, \;C_r^{KM} QL\colon=  s(C_r(s.) QL)_{s.\in S_K(M)}  
\leqno{(16)}$$
 In the case  $ K = M $ we write  
$C_r^{K} L $ (resp. $C_r^{K} ICL, \; C_r^{K} QL $ for $ C_r^{KK} L $ (resp. $ICL, \; QL $ 
instead of $L$ }
 $$C_r^K L = s(C_r(J,s.))_{(J,s.) \in K_.^+ \times S(K)} = s(
(\cap_{K\supsetneqq s_{\lambda}\in s.} 
W^{s_{\lambda}}_{a_{s_{\lambda}}(J,r-1)}Gr^{W^K}_{a_{K(J,r)}}L)_{(J,s.)
\in K_.^+ \times S(K)} \leqno{(17)}$$
\section {Main theorems on the properties of the weight filtration} In
this section we aim to prove that the filtration is actually the weight of what would be 
in the proper case
 a mixed Hodge complex in Deligne's terminology, that is the induced filtration by $F$ on 
the graded parts $Gr^{\cal W}$ of ${\cal W}$ is a Hodge filtration. For this we need to 
decompose the complex as a direct sum of intermediate extensions of variations of Hodge 
structures (which has a meaning locally) whose cohomologies are  pure Hodge
structures [5] and [24] in the proper case. This is done in the following three 
subsections. In the first we prove a  key lemma that  apply to prove the purity of the 
complex $C^K_r L$. Once this purity result is established, we can easily prove in the 
second subsection the decomposition theorem after a careful study
of the category of indices $S(I)$. In the third subsection we give the global statements 
for a filtered combinatorial logarithmic complex. For this we use the above local 
decomposition to obtain a global decomposition 
of the graded weight into intermediate extensions of polarised $VHS$ on the 
various intersections of components of $Y$. This last statement uses the formula 
announced by Kashiwara and Kawai [26] that we prove since we have no reference for its 
proof.

\subsection  {Purity of the cohomology of the complex $C^K_r L$} In this
subsection we introduce  the fiber of  the variations of Hodge
structures  needed in the decomposition of $Gr^{\cal W}$. The result here is the 
fundamental step in the general proof. The plan of this subsection is as follows. First 
we start with a key lemma relating the
various relative monodromy weight fltrations (centered at zero) associated to a nilpotent 
orbit $L$; $N_i$ is compatible with $W(N_j)$
for $i \neq j$ but shift Hodge filtration by $-1$, hence it is not clear whether it is
strict, however we need technical results to establish the purity and decomposition  
properties  and this key lemma provides what seems to be
the elementary property at the level of a nilpotent orbit that leads to the 
decomposition. Second we present a set of elementary complexes. Finally we
state the purity result on the complexes $C^K_r L$ which behave as a direct sum of 
elementary complexes.
\be{Properties of the relative weight filtrations}{2.I.1}  \

Given the nilpotent orbit we may consider various filtrations $W^J = W(\Sigma_{i\in 
J}N_i)$ for various $J \subset M$. They are  centered at $0$, preserved by 
$N_i$ for $i\in M$ and shifted by $-2$ for $i \in J$: $N_iW^J_r \subset
N_iW^J_{r-2}$. We will need more subtle relations between these
filtrations that we discuss in this subsection.

\smallskip

{\em Key lemma (Decomposition of the relative weight filtrations) : Let
$ (L, N_i, i \in [1,n], F)$ be a polarised nilpotent orbit 
and let $W^i \colon= W(N_i)$ (all weights centered at $0$ ), then :

\noindent i) For each subset $ A = \{ i_1,\ldots, i_j\}\subset [1,n]$,
of length  $\mid A \mid = j $ $$ Gr_r^{W^A} L \simeq \oplus_{m.\in
X^A_r} Gr^{W^{i_j}}_{m_{i_j}} \cdots 
Gr^{W^{i_1}}_{m_{i_1}} L \;{\hbox {where}}\; X^A_r = \{m.\in {\Z}^j: 
\Sigma_{i_l \in A } m_{i_l} = r\} $$   
and if $A = B \cup C $
$$ Gr_{b+c}^{W^A} Gr_c^{W^C} L \simeq Gr_b^{W^B} Gr_c^{W^C} L \simeq
Gr_{b+c}^{W^A} 
Gr_b^{W^B} L$$ 
\noindent  ii) Let $ N_i' $ denotes the restriction of $N_i$ to $Gr^{W^C}_c$ and 
$ N_B' = \Sigma_{i\in B}N'_i$, then  $W_b^B $ induces $W_b( N_B')$ on $Gr^{W^C}_c$.

\noindent iii) In particular the repeated graded objects in (i) do not depend on the 
order of the elements in $A$.}

\smallskip

\noindent {\em Remark: This result give relations between  various 
weight filtrations in terms of the elementary ones $W^i \colon= W(N_i)$
and will be extremely useful in the study of the properties of the weight filtration on 
the higher direct image.
For example for $N = N_1 +  \cdots + 
N_n$, there exists a canonical decomoposition }
 $$ Gr_r^{W(N)}L = \oplus_{m.\in X_r} Gr^{W^n}_{m_n} \cdots
Gr^{W^1}_{m_1} L 
\;{\hbox {where}}\; X_r = \{m.\in {\Z}^n: \Sigma_{i \in [1,n]} m_i = r \}.   $$   
The proof by induction on $n$ is based on the following important result of Kashiwara 
[25, thm 3.2.9, p 1002]:

\smallskip

\noindent {\em Let $(L,N, W) $ consists of a vector space endowed with
an increasing filtration $W$ preserved by a nilpotent  endomorphism  $N$
on $L$ and suppose that the relative filtration $ M = M(N,W)$ exists, 
then there exists a canonical decomposition:}  $$Gr^M_l L = \oplus_k Gr^M_l Gr^W_k L .$$ 
{\it Proof of the key lemma}.  To stress the properties of commutativity of the graded 
operation for 
the filtrations, we prove first 

\smallskip

\noindent {\em Sublemma:  For all subsets $[1,n]\supset A \supset
\{B,C\},$ the 
isomorphism of Zassenhaus
$Gr^{W^B}_b Gr^{W^C}_c L \simeq Gr^{W^C}_c Gr^{W^B}_b L $
is an isomorphism of $MHS$ with weight filtration (up to a shift) $W = W^A$ and Hodge
filtration $F =  F_A $, hence compatible with the third filtration $W^A$
or $F_A$.}

\smallskip

Proof: Recall that both spaces $Gr^{W^B}_b Gr^{W^C}_c $ and $ Gr^{W^C}_c
Gr^{W^B}_b$ are isomorphic to $W^B_b \cap W^C_c $ modulo $ W^C_c \cap W^B_{b-1} + W^B_b 
\cap W^C_{c-1} $. In this isomorphism a third filtration like
$F_A$ (resp.
$ W^A $) is induced on one side by $ F'_k = (F_A^k \cap W^C_c) + W^C_{c-1}$ (resp.  $W'_k 
= W_k^A \cap W^C_c) +  W^C_{c-1}$ ) and on the 
second side by  $ F_k'' = (F_A^k \cap W^B_b) +  W^B_{b-1}$ (resp. $W_k'' = W_k^A \cap 
W^B_b) +  W^B_{b-1}$ ). We 
introduce the third  filtration $ F_k''' = F_A^k \cap W^B_b \cap W^C_c $ (resp. $W_k''' = 
W_k^A \cap W^B_b \cap W^C_c $) 
and we notice that all these spaces are in the category of $MHS$, hence the isomorphism 
of Zassenhaus
 which must be strict, is compatible with the third 
filtrations induced by $F_A$ ( resp. $W^A$ ).

\smallskip

Proof. i) Let $ A \subset [1, n]$ and $i\in A$, then $W^A $ exists and induces the 
relative weight filtration for $N_i$ with respect to $W^{(A-i) }$. Then we have by 
Kashiwara's result
$ Gr_l^{W^A}L = \oplus_k  Gr_l^{W^A} Gr_k^{W^{(A-i)}}L  $. Let us attach to each point 
$(k,l)$ in the plane the space $ Gr_l^{W^A} Gr_k^{W^{(A-i)}}L  $ and let  
$M_j = \oplus_l  Gr_l^{W^A} Gr_{l-j}^{W^{(A-i)}}L $ be the direct sum along 
indices on a parallel to the diagonal (shifted by $j$) in the 
plane $(k,l)$. Then we have for $j>0$
$$ (N_i)^j: Gr_{k +j}^{W^A} Gr_k^{W^{(A-i)}}L \simeq    Gr_{k-j}^{W^A} Gr_k^{W^{(A-i)}}L, 
\; (N_i)^j: M_j \simeq M_{-j}.$$
This property leads us to introduce 
the space  $V = \oplus_l  Gr_l^{W^A } L \simeq \oplus_{l,k}  Gr_l^{W^A} Gr_k^{W^{(A-i)}}L 
$, then $N_i$ on $ L$ extends to a nilpotent endomorphism on $V$,  $ N_i: V \rightarrow V 
$ inducing $ N_i:Gr_l^{W^A}L 
\rightarrow Gr_{l-2}^{W^A }L $ on each $l-$component of $V$. We consider on $V$ two 
increasing filtrations  $W'_s \colon= \oplus_{l-k \leq s}
Gr_l^{ W^A }Gr_k^{W^{(A-i)}}L$ and 
  $ W_s'' \colon= \oplus_l W_s^i Gr_l^{W^A}L $. Then
$ N_i$ shift these filtrations by $-2$. In fact $N_i: W'_s \rightarrow
W'_{s-2}$ sends $Gr^{W^A}_lGr_k^{W^{(A-i)}}L$ to
$Gr^{W^A}_{l-2}Gr_k^{W^{(A-i)}}L$ and $ (N_i)^j $ induces an isomorphism
$ Gr^{W'}_j V \simeq Gr^{W'}_{-j} V$. As well we have an isomorphism $
Gr^{W''}_j V \simeq Gr^{W''}_{-j} V$, since 
$(N_i)^j : (Gr^{W^i}_j L, W^A, F_A) \simeq  (Gr^{W^i}_{-j} L, W^A, F_A)
$ is an isomorphism of $MHS$ up to a shift in indices, hence strict on $
W^A $ and $ F_A$ and induces an isomorphism  $ Gr^{W^i}_j Gr_l^{W^A} 
\simeq Gr^{W^i}_{-j}Gr_{l-2j}^{W^A} $.
  Then these two filtrations $W'_s $ and $W_s''$ are equal by uniqueness
of the weight filtration of  $ N_i$ on $V$, that is 
$$ W_s'' = \oplus_{k,l} W^i_s Gr_l^{W^A } Gr_k^{W^{(A-i)}}L = W'_s =
\oplus_{l-k \leq s} Gr_l^{W^A} Gr_k^{W^{(A-i)}}L $$ that is $ W^i_s
Gr_l^{W^A } Gr_k^{W^{(A-i)}}L =  Gr_l^{W^A } Gr_k^{W^{(A-i)}}L $ if $l-k
\leq s$ and $ W^i_s Gr_l^{W^A } Gr_k^{W^{(A-i)}}L = 0 $ otherwise, or $$
Gr_l^{W^A}  Gr^{W^i}_j L \simeq  Gr_l^{W^A} Gr_j^{W^i} 
Gr_{l-j}^{W^{A-i}}L, {\hbox {and for all}} \; l \neq k + k', \;
Gr_l^{W^A} Gr_k^{W^i} Gr_{k'}^{W^{(A-i)}}L \simeq 0. $$ In other words:
$W^A$ induces a trivial filtration on $Gr_k^{W^i} Gr_{k'}^{W^{(A-i)}}L $
of weight $k+k'$ that is 
$$ Gr_l^{W^A} L \simeq \oplus_k Gr_l^{W^A} Gr_k^{W^{A-i}} L  \simeq \oplus_k Gr_l^{W^A} 
Gr_k^{W^{A-i}} Gr^{W^i}_{l-k} L \simeq \oplus_k Gr_k^{W^{(A-i)}} Gr_{l-k}^{W^i} L.$$ 
Now if we suppose by induction on length of $A $, the decomposition true for $ A-i$, we 
deduce easily the decomposition for $ A$ from the above result.

\noindent ii) We restate here the property of the relative monodromy for $W^A$ with 
respect to $W^C$.
 
\noindent iii)In the proof above we can start with any $i$ in $A$, hence the 
decomposition 
is symmetric in elements in $A$. It follows that the graded objects of the 
filtrations $W^i, W^r, W^{\{i,r,j\}}$ commute 
and since $W^j$ can be expressed using these filtrations, we deduce that
$W^i, W^r, W^j$ also commute, for example: $ Gr_{a+b+c}^{W^{\{i,r,j\}}} 
Gr_{a+b}^{W^{\{i,j\}}} Gr_a^{W^i} \simeq Gr_c^{W^r}Gr_b^{W^j}
Gr_a^{W^i}$ is symmetric in $i,j,r$.

\smallskip
\noindent {\em Corollary:
The morphism $N_i$ induces for all $j$, exact sequences for all integers $r$ :
 
$0\rightarrow W_r^j\cap ker N_i \rightarrow W_r^j L \rightarrow
W_r^j\cap N_i L \rightarrow 0 $ 
and
$0\rightarrow W_r^j \cap  N_i L \rightarrow W_r^j L \rightarrow W_r^j (L/N_i L 
)\rightarrow 0 .$}

\smallskip
Proof:
$N_i$ is strict for $W^i $ and $ W^{\{i,j\}}$ hence we have the above exact sequences 
for $ Gr^{  W^{\{i,j\}}}_{a+b} Gr^{W^i}_b =  Gr^{W^{j}}_a Gr^{W^i}_b $.

\be {Elementary complexes}{2.I.2}
The proof of the purity uses the following simplicial vector spaces. For
each $J \subset [1,n], $ let
$$ K((m_1,\cdots,m_n),J) L = Gr^{W^n}_{m_n - 2 \mid n \cap J \mid }
\cdots Gr^{W^r}_{m_r - 2 \mid r 
\cap J \mid } \cdots Gr^{W^1}_{m_1 - 2 \mid 1 \cap J \mid } L$$
 (resp. for $L/N_J L $ and $ N_J L $). For all $i \notin J$, the endomorphism $ N_i $ 
induces a morphism 
denoted also $ N_i: K((m_1,\cdots,m_n),J) L \rightarrow
K((m_1,\cdots,m_n),J\cup i)L,$  
(resp. for $L/N_J L $ and $ N_J L $ instead of $L$),
then we consider the following elementary complexes defined as simple associated 
complexes:
$$K(m_1,\cdots,m_n) L \, \colon= s( K((m_1,\cdots,m_n),J) L , N_i )_{J \subset 
[1,n]}\leqno{(18)}$$ 
(resp. $K(m_1,\cdots,m_n) QL \, \colon= s(
K((m_1,\cdots,m_n),J) L/N_J L , N_i )_{J  \subset [1,n]}$ 

\noindent and $K(m_1,\cdots,m_n) IC L \, \colon= s(
K((m_1,\cdots,m_n),J) N_J L  , N_i )_{J  \subset [1,n]}).$ 

\smallskip 
\noindent {\em Proposition: i) For any $((m_1,\cdots,m_n) \in \Z^n $ let $J(m.) = 
\{i \in  [1,n]: m_i \geq 1\}  $, then  the cohomology of an elementary complex
 $K(m_1,\cdots,m_n) L$ is a subquotient of $K((m_1,\cdots,m_n), J(m.))L$, 
hence concentrated in degree $\mid J(m.) \mid $. Moreover it vanishes
iff there exists at least one $ m_i = 1$. More precisely, if no $ m_i = 1$, the 
cohomology is isomorphic to
  $ \;\;  K((m_1,\cdots,m_n),J(m.) ) [(\cap_{i \notin J(m.)} (ker N_i:
L/(\Sigma_{j \in J(m.)  } N_j L) \rightarrow L/(\Sigma_{j \in J(m.)} N_j L))]$, 

\noindent moreover this object is symmetric in the operations kernel and
cokernel and is isomorphic to 

\noindent $K((m_1,\cdots,m_n),J(m.))  [(\cap_{\{i: m_i < 1 \}} ker 
N_i)/(\Sigma_{\{j : m_j > 1 \}} N_j (\cap_{\{i: m_i < 1 \}} ker N_i))]$,

\noindent that is at each process of taking $Gr_{m_i}^{W^i} $ we apply the functor $ ker 
$ if $m_i \notin J(m.) $ and $coker $ if $m_i \in J(m.)$.

\noindent ii) If there exists $m_i > 0$, then $ K(m_1,\cdots,m_n) IC L
\simeq 0,$ hence 

\noindent $ K(m_1,\cdots,m_n)  L \simeq K(m_1,\cdots,m_n) Q L .$

\noindent iii) If all $m_i \leq 0 $, then 

\noindent $K(m_1,\cdots,m_n) IC L \simeq
K(m_1,\cdots,m_n)  L \simeq Gr_{m_n}^{W^n} \cdots Gr_{m_1}^{W^1} \cap_{i
\in [1,n]} ker ( N_i: L \rightarrow L )$.}

\smallskip

 \noindent Proof. i) It is enough to notice that $N_i$ is injective if $m_i > 0$, 
surjective if $m_i < 2$ and bijective if $m_i =1$.
In fact, given $N_i$ we can view $K(m_1,\cdots,m_n) L$ as the cone over

 \noindent $N_i:
K(m_1,\cdots,\widehat{m_{ i}},\cdots, m_n) (Gr^{W^i}_{m_i}L)\; {\buildrel N_i
\over 
\rightarrow } \;K(m_1,\cdots,\widehat{m_{i}},\cdots, m_n) (Gr^{W^i}_{m_i-2}L)$. 

\noindent Hence if $m_i > 0 $ (resp. $m_i < 2$),
that is $ N_i $ is injective on $Gr_{m_i}^{W^i} $ (resp. surjective onto 
$Gr_{m_i-2}^{W^i} $ ),

\noindent $ K(m_1,\cdots,m_{i},\cdots, m_n) L
\cong  K(m_1,\cdots,\widehat{ m_{i}},\cdots, m_n) (Gr^{W^i}_{m_i - 2}(L/N_iL))[-1]$

\noindent (resp. $K(m_1,\cdots,\widehat{m_{i}},\cdots, m_n) (Gr^{W^i}_{m_i}(ker N_i: 
L\rightarrow  L)$)

\noindent where $K(m_1,\cdots,\widehat{m_{i}},\cdots, m_n) $ is viewed for the nilpotent 
orbit $Gr^{W^i}_{m_i-2} L/N_i L$ 

\noindent (resp.$ Gr^{W^i}_{m_i}(ker N_i: L\rightarrow  L)$) with the nilpotent 
endomorphisms $ N_j'$ induced by $N_j$ for $j \neq i$.

\noindent ii) $ K(m_1,\cdots,m_n) IC L $ can be viewed as a cone over 

$ N_i : K(m_1,\cdots,\widehat{m_{i}},\cdots, m_n) IC(Gr^{W^i}_{m_i} L)
\rightarrow 
K(m_1,\cdots,\widehat{m_{i}},\cdots, m_n) IC (Gr^{W^i}_{m_i - 2} N_i L) $ 

\noindent where
$K(m_1,\cdots,\widehat{m_{i}},\cdots, m_n) (Gr^{W^i}_{m_i}L$ is viewed for the
$(n-1)-dim$ nilpotent orbit $ Gr^{W^i}_{m_i} L$. If $ m_i  > 0 $, then $ N_i $ is an 
isomorphism.

\noindent iii) If $m_i \leq 0$, then the above $N_i$ is just surjective and 
$ K(m_1,\cdots,m_n) IC L $ is isomorphic to 
 $ K(m_1,\cdots,\widehat{m_{i}},\cdots, m_n) IC(Gr^{W^i}_{m_i} ker ( N_i : L
\rightarrow N_i L )$, then (iii) follows by induction on $i$.

\be {Main result}{2.I.3} \
\noindent {\em Theorem (Purity). Let $L$ be a polarised nilpotent orbit ( local 
hypothesis 2 
(\S 1, I.3)), then the complexes $C^K_r L $ in (\S 1, II.3) where we suppose $K = M $ of 
length $\mid K \mid = n $, satisfy the following properties

 i) Let $r > 0$ and
$T(r) = \{(m_1,\cdots,m_n) \in {\bf N}^n: \forall i \in K, m_i \geq 2 $
and $ \Sigma_{j \in K } m_j  = \mid K \mid +r  \}$

 \noindent then the cohomology of the complex $C_r^K L$ is isomorphic to that of the 
following complex
$$C_r^K L \cong C(T(r))\simeq  \oplus_{(m_1,\cdots,m_n) \in T(r)} K(m_1, \ldots, m_n) L$$ 
\noindent In particular its cohomology, concentrated in degree $\mid K \mid $, is 
isomorphic to 

 $ Gr^{W^K}_{r- \mid K \mid}[L/(\Sigma_{i \in K} N_i
L)] \simeq
\oplus_{(m_1,\cdots,m_n) \in T(r)} Gr^{W^n}_{m_n - 2 } \cdots 
Gr^{W^i}_{m_i - 2 } \cdots Gr^{W^1}_{m_1 - 2 } [L/(\Sigma_{i \in K} N_i
L)] $

\noindent it is a polarised Hodge structure of  weight $r + m - \mid K \mid$ with the 
induced filtrations $W^K$ (shifted by $m$) and $F^K$.
Moreover, if $r = 0$, the complex $C_r^K L$ is acyclic. 

ii) Dually, for $ r <0 $, let
 $T'(r) = \{(m_1,\cdots,m_n) \in {\Z}^n: \forall i \in K, m_i \leq 0, \Sigma_{j \in K } 
m_j  = \mid K \mid +r  \},$
 then the complex $ C_r^K L$ is isomorphic to the  following
complex
$$C_r^K L \cong C(T'(r))[1 - \mid K \mid] \simeq
\oplus_{(m_1,\cdots,m_n) \in T'(r)} K(m_1, \ldots,
m_n) L [1 - \mid K \mid]$$
  In particular its cohomology, concentrated in degree $\mid K \mid - 1$, is isomorphic 
to 

\noindent $ Gr^{W^K}_{r+\mid K \mid}[(\cap_{i \in K} (ker
N_i: L \rightarrow L ] \simeq
\oplus_{(m_1,\cdots,m_n) \in T'(r)} Gr^{W^n}_{m_n } \cdots 
Gr^{W^i}_{m_i } \cdots Gr^{W^1}_{m_1 } [(\cap_{i \in K} (ker N_i: L \rightarrow L ] $

\noindent it is a polarised Hodge structure of weight $r + m + \mid K \mid$  
with
the 
induced filtrations $W^K$ (shifted by $m$) and $F^K$ .

iii) The complex $C_r^K L$ is quasi-isomorphic to $C_r^K Q $ (16) for $r \geq 0 $ 
and to the complex $C_r^K IC L$  (16) for $ r \leq 0.$} 

\smallskip

\noindent {\em Remark:  If $ r \in [1, \mid K \mid - 1]$, $T(r)$ is empty and $C^K_r L$ 
is acyclic.  If $ r \in [ - \mid K \mid + 1, 0]$, $T'(r)$ is empty and $C^K_r L$
is acyclic. In all cases $C^K_r L$ appears in $Gr^{{\cal W}}
{\Omega}^*({\cal L})$.}

\smallskip
\noindent Proof:
 The important fact used here is the particular decomposition  for a
nilpotent 
orbit of the relative filtrations, that is the isomorphism, functorial
for the 
differentials of $C^K_r L$
$$ Gr^{W^K}_{a_{K(J,r)}} ( \cap _{K\supsetneqq s_{\lambda}\in s.} 
W^{s_{\lambda}}_{a_{\lambda}(J,r-1)}) L \simeq 
 \oplus_{m. \in X(J,s.,r)} Gr^{W^n}_{m_n}\cdots Gr^{W^i}_{m_i  } \cdots Gr^{W^1}_{m_1} 
L,\;   $$
where for all $ (J,s.)\in K_.^+ \times S(K), $ 
$$ X(J,s.,r) = \{ m. \in {\Z}^n :\Sigma_{i \in K} m_i = a_{K}(J,r)\; \hbox {and} \;
\forall s_{\lambda}\in s., \Sigma_{i \in s_{\lambda } } m_i \leq a_{\lambda}(J,r-1)  \}$$

\noindent In particular, if we define for $J = \emptyset, 
 X(s.,r) =  X(\emptyset,s.,r) $, 

\noindent $ X(s.,r) = \{ m. \in {\Z}^n :\Sigma_{i \in K} m_i = \mid K \mid + r \;\; \hbox 
{and} \;\;
\forall s_{\lambda}\in s., \Sigma_{i \in s_{\lambda } } m_i \leq \mid s_{\lambda}\mid  + 
r-1  \}$,

\noindent  the complex $C_r^K (s.)$ splits as a direct sum of elementary
complexes 
$$ C_r^K (s.)\simeq \oplus_{ X(s.,r)}  K(m_1,\ldots,m_n).$$ \noindent On
the otherside for a fixed $J \subset K $ we consider the complex defined by the column of 
vector spaces
 $$ C_r^K (J) = s[C_r^K(J,s.)]_{s. \in S(K)}\cong  s[ \oplus_{m. \in X(J,s.,r)} 
K((m_1,\ldots,m_n),J)]_{s.\in S(K)}.$$
We want to show that each column is an acyclic complex if
$((m_1,\cdots,m_n) \notin T(r)$ and a resolution of
$K((m_1,\ldots,m_n),J)$ otherwise (if
$((m_1,\cdots,m_n) \in T(r)$). 

\noindent This is just a combinatorial study, which helped to formulate the statement 
after an explicit study of the theorem in case $n =2$ and $n = 3 $. We give a proof based 
on the following facts:

\medskip
\noindent {\em Lemma: i) Let $r \geq 0$, then for each $i\in K$ the
sub-complexes
$$C_r^K (W_{1}^iL) \simeq s(\oplus_{m.\in  X(s.,r)\;\hbox {and}\;
 m_i < 2}  K(m_1,\cdots,m_n) L )_{s.\in S(K)} $$

\noindent as well $C_r^K (IC (W_{1}^i L )) $ are acyclic.
More precisely for $r > 0 $, each column 

\noindent $C_r^K (J) (W_{1}^iL) =
s(C_r^K(J,s.) (W_{1}^iL)_{s. \in S(K)} $  is acyclic. 

ii) Dually, for $r < 0$  and for each $i\in K$ the quotient complexes
$$C_r^K (L/W_{1}^iL) \simeq s(\oplus_{m.\in  X(s.,r)\;\hbox {and}\; m_i
\geq 2}  K(m_1,\cdots,m_n) L )_{s.\in S(K)} $$

\noindent are acyclic column by column.}
 
\noindent Proof: We  distinguish in $S (K)$ the subcategory $S''_i$
whose objects $s.$ contain $K$ and $K-i$. The complement $ S(K) - S''_i
$ is a full subcategory of $S(K)$ since if we delete a subset in $s. \in
S(K) - S''_i $  we still have an object in this subcategory. Hence the
sum $C''_i \colon=
s[C_r^K(s.)]_{s. \in S(K) - S''_i} $   is a subcomplex of $C^K_r L $
whose quotient complex is 
$C(S''_i) \colon=  s[C_r^K(s.)]_{s. \in S''_i}. $ 

Dually, we consider  $S'''_i \subset S'_i$  whose objects $s.$ contain
$\{i\}$. The complement $ S(K) - S'''_i $ is a full subcategory of
$S(K)$ since if we delete a subset in $s. \in S(K) - S'''_i $  we still
have an object in this subcategory. Hence the sum $C'''_i \colon=
s[C_r^K(s.)]_{s. \in S(K) - S'''_i} $   is a subcomplex of $C^K_r L $
whose quotient complex is 
$C(S'''_i) \colon=  s[C_r^K(s.)]_{s. \in S'''_i}. $    

\medskip
\noindent {\em Sublemma: In the exact sequence
$$0\rightarrow C''_i(W_1^iL) \rightarrow C_r^K(W_1^iL)[1] \rightarrow
C(S''_i)(W_1^iL)[1] \rightarrow 0 $$
 the complexes at each side are acyclic (column by column if $r>0$), so
is the middle complex.
Dually, in the exact sequence
$$0\rightarrow C'''_i(L/W_0^iL) \rightarrow C_r^K(L/W_0^iL)[1]
\rightarrow 
C(S'''_i)(L/W_0^iL)[1] \rightarrow 0 $$
 the complexes at each side are acyclic (column by column if $r < 0$),
so is the 
middle complex.}

Proof. We write $s.$ as $(s.'\supset s_v \cup i \supset s_{v-1}\supset 
s.'')$ where$ i \notin s_{v-1} $ and   distinguish in the objects of $S
(K)$ two families : $S_i $ whose objects are defined by the $s.$
satisfying 
$ s_v \supsetneqq s_{v-1}$ (including the case  
$s_{v-1} = \emptyset $)  and $S'_i$ whose objects $s.$ satisfy $s_{v-1}
= s_v $ 
(including the case  $s_{v-1} = \emptyset $, that is $s_v \cup i = i) $.

We form the complexes $C(S'_i - S''_i) \colon=  s[C_r^K(s.)]_{s. \in
S'_i - S''_i} $ (resp. $C(S_i)\colon=  s[C_r^K(s.)]_{s. \in S_i }  $. If
we delete $s_v \cup i $ in $s. \in S'_i - S''_i$ we get an element in
$S_i$ . then removing $s_v \cup i $ can be  viewed as a morphism $
d_{s_v \cup i
}: 
C(S'_i - S''_i)  \rightarrow C(S_i) $ and the cone over this $ 
d_{s_v \cup i }$ is equal to $C''_i[1] $.

Now, if we reduce the construction to $W_1^iL$  and if
 $s.$ is an object of $S_i'$, the condition $m.\in  X(s.,r)$ associated
to $s.$ when $ s_{v-1} = s_v \neq  \emptyset$ (resp. $s_v \cup i$) is
$\Sigma_{j \in s_v} m_j \leq  \mid s_v \mid + r-1 $ (resp. $m_i +
\Sigma_{j \in s_v} m_j \leq  \mid s_v \cup i \mid + r-1 $), but
precisely when $m_i < 2$ (that is in $W^i_1 L $ the condition for $s_v$
is equivalent to the union of the conditions for $s_v$ and $s_v \cup i$,
that is 
$ d_{s_v \cup i } $ 
induces an isomorphism for such object in $S'_i - S''_i$.

When $s_{v-1} = \emptyset$ and if $r>0$, the condition $m_i \leq r$ is
irrelevant since already 
$ m_i \leq 1 $ and $ r \geq 1 $, so that  $ d_{s_v \cup i } $ 
induces an isomorphism for all objects in $S'_i- S''_i$ (if $r =0$ the
difference 
are 
complexes $K(m_1, \ldots, m_n)$ with some $m_i = 1 $, hence acyclic ).

Dually, $S'''_i $ whose objects $s.$ satisfy
$s_{v-1} = s_v  = \emptyset $, that is $ i \in s.$ is contained in $
S'_i$ , then the cone over the morphism $ d_{s_{v-1} }: C(S'_i - 
S'''_i)(L) \rightarrow C(S_i)(L)$ is isomorphic to $C'''_i [1].$ 
When $m_i > 0$ the condition for $s_v \cup i$ is equivalent to the union
of the conditions for $s_v = s_{v-1} \neq  \emptyset$  and $s_v \cup i$,
that is $ d_{s_{v-1 }} $ induces an isomorphism for such object in
$S'_i-S'''_i$.

Finally we prove $C(S''_i)(W_1^iL) = 0.$ In fact the conditions for $K$
and 
$K-i$
in any $s. \in S''_i$ are 
$\Sigma_{j \in K} m_j = \mid K \mid + r $
(resp. $ \Sigma_{j \in ( K-i )} m_j \leq  \mid K-i  \mid + r-1 $), hence
we get 
by difference $m_i \geq 2$.

Dually,$C(S'''_i)(L/W_0^iL) = 0.$ In fact the condition  $ m_i \leq r$
corresponding to $i \in s. \in S'''_i $ is not compatible with $m_i > 0$
when $ r < 0 $ which ends the proof of the sublemma.

On the otherside, it is easy to check that 

\smallskip

\noindent {\em Lemma: The complex $C(T(r))$ 
is contained in each $C^K_r(s.)$ that is $T(r)\subset  X(s.,r)$.
Dually, the complex $C(T'(r))$ 
is contained only in $C^K_r(s.)$ for $s. = K$.}

\smallskip

\noindent  We check the condition 
$\forall s_{\lambda} \in s., \Sigma_{i \in s_{\lambda } } m_i \leq \mid 
s_{\lambda}\mid + r-1  \}$ for all $m. \in T(r) $ by induction: 
$\Sigma_{j \in K } m_j  = \mid K \mid +r \Rightarrow \forall A = K - k \subset K, 
\Sigma_{j \in A } m_j \leq \mid A\mid  + r-1  $ by substracting   $m_k > 1$.
If this is true for all $A : \mid  A \mid  = a $ then $ \forall B = A - k\subset 
A, \Sigma_{j \in B } m_j \leq \mid B \mid  + r- 2  $ as well.

Dually, the condition for $K$, $\Sigma_{j \in K } m_j  = \mid K \mid +r 
\Rightarrow \forall A = K - k \subset K, 
\Sigma_{j \in A } m_j > \mid A\mid  + r-1  $ by substracting   $m_k < 1$.
If this is true for all $A : \mid  A \mid  = a $ then $ \forall B = A - k\subset 
A, \Sigma_{j \in B } m_j > \mid B \mid  + r- 1  $ as well.

Finally, we form the complex $C(r) = s(C(T(r))_{s. \in S(K)} \subset C^K_r L $ and we 
deduce

\smallskip

\noindent {\em Lemma: i) The quotient $(C^K_r L / C(r)) \cong 0 $ is
acyclic.

ii) For any maximal index $s.\in S(K)$, the embeddding $C(T(r)) \subset C^K_r(s.)$ 
induces a quasi-isomorphism $C(T(r))\cong C^K_rL $.

iii) Dually, the quotient $[C^K_r L / C(T'(r))(s. = K)] \cong 0 $ is acyclic. }

\smallskip

Proof. i) Since any element of the quotient can be represented by an
element in 
a subcomplex $K(m_1,\ldots, m_n)$ with some $m_i < 2$ we can apply the
lemma  above for some $W^i_1 L$.

ii) $ d_{s_v \cup i }$ induces an isomorphism on the complexes obtained
as sum of $C(T(r))$ over $S'_i - S''_i$ and $S_i$, hence the cohomology 
of $C(r)$ comes from $C(S''_i)$. All elements $s.$ in $S''_i$ contain $K
\supset K-i$, so we can repeat the same arguments  in the category
$S(K-i) $ but 
for $j \neq i$, so $(ii)$ follows by induction.

iii) The assertions for $C_r^K IC L $ and $ C_r^K Q L$  (20) follow from the
the corresponding isomorphisms $ C(T(r))L \simeq (C(T(r))IC L $ and 
$(C(T'(r))L \simeq (C(T'(r))IC L $. 

\medskip

{\em Remark (duality)}. Given a polarised nilpotent orbit $(L, N_i (
i\in M), P)$, the local duality induces an isomorphism:

$$d^K_r : C^K_r L [2 \mid K\mid - 1]  \simeq  Hom_{\Q}(C^K_{-r} L, \Q)
$$ hence $H^i (C^K_r L)[2 \mid K\mid - 1]) = H^{i + 2 \mid K\mid - 1}
\simeq  Hom_{\Q}(H^{-i}(C^K_{-r} L), \Q)  $ and for $i = 1 - \mid K \mid
$: 
$$d^K_r : H^{\mid K\mid } (C^K_r L) \simeq  Hom_{\Q}(H^{\mid K\mid -1}
(C^K_{-r} L), \Q) $$

The duality is constructed as follows:
For each $s.$ we define $C(s.) = \{ s'. \in S(K): \mid s'.  \mid + \mid
s.\mid =\mid K \mid + 1$ and $ s. \cup s'.$ is maximal $ \}$, that is to
say $ s.$ is complementary to $s'.$ except that both must contain $K$,
then we define:
$P^* : C_r(J,s.)\otimes C_{-r}(J',s'.))\rightarrow \Q $ 
as $ P^* (a, b) = P (a,b)$ for $J' = K-J $ and $s'. \in  C(s.)$ and zero
otherwise 
($P^* $ is non zero on $ C_r(J,s.)\otimes (\oplus_{s'. \in C(s.)
}C_{-r}(K-J,s'.)) $. 

It can be checked that the induced morphism $P^*: C^K_r L [2 \mid K\mid
- 1] \rightarrow Hom_{\Q}(C^K_{-r} L, \Q)  $ commutes with
differentials.

\smallskip 

\noindent {\em The relation between $ C^K_r L$ and $C^{KM}_r L$} 

\smallskip

The following result will be important in  the general proof
of the decomposition of $Gr^{\cal W}_r  \Omega^* L $ as direct sum of intersection 
complexes.

\smallskip

\noindent {\em Proposition: 
Let $H = H^* (C_r^K L$), considered as a nilpotent orbit with indices $i \in M-K$, then 

\noindent i) We have : $C^{KM}_r L \simeq {\cal W}_{-1 }\Omega^*(H)$

\noindent ii) For $r \geq 0 $, $ C^{KM}_r L \simeq C^{KM}_r QL $}

\smallskip

\noindent Proof. i) We can write $(K_.^+ \times S(K) \times (M-K)_.^+
\times S(M-K)) \simeq M_.^+\times S_K(M)$, with the correspondence $ 
(J,s.,J',s.')\rightarrow ((J,J'),(s'_.\cup K, s.)) $  then using the relations:

\noindent 1)$a_{s_{\lambda}}((J,J'),r-1) = a_{s_{\lambda}}(J,r-1)$ 
when $ s_{\lambda}\subset K $ 

\noindent 2)$a_{s'_{\lambda}}((J,J'),r-1) = a_{s'_{\lambda}}(J',-1) + a_K(J,r)$
and 

\noindent $W^{s'_{\lambda}\cup K}_{a_{s'_{\lambda}\cup K}((J,J'),r-1)} 
(Gr^{W^K}_{a_K(J,r)}) = W^{s'_{\lambda}}_{a_{s'_{\lambda}}(J',-1)} (Gr^{W^K}_{a_K(J,r)}) 
$ since $W^{s'_{\lambda}\cup K}$ is relative to
$W^K$

we find
 
 \noindent $C_r^{KM}L \colon=  s( \cap _{K\neq s_{\lambda}\in (s.'\cup 
K, s.)}W^{s_{\lambda}}_{a_{s_{\lambda}}((J,J'),r-1)}
Gr^{W^K}_{a_K((J,J'),r)} 
L)_{((J,J'),(s'_.\cup K, s.)) \in M_.^+\times S_KM} \simeq $

\noindent $s[ (\cap _{ s'_{\lambda}\in 
s.'}W^{K \cup s'_{\lambda}}_{a_{s'_{\lambda}}(J', a_K(J, r-1)} (s [
Gr^{W^K}_{a_K(J,r)}(\cap _{K\supsetneqq s_{\lambda}\in s.} 
W^{s_{\lambda}}_{a_{s_{\lambda}}(J,r-1)} L ]_{(J,s.)\in 
K_.^+ \times S(K)}) ]_{(J',s.')\in (M-K)_.^+ \times S(M-K)}$ 

\noindent $ \simeq s[ \cap _{ s'_{\lambda}\in s.'}
W^{s'_{\lambda}}_{a_{s'_{\lambda}}(J',-1)}  ( C_r^K L)]_{(J',s.')\in (M-K)_.^+ \times 
S(M-K)} $

\noindent where  $\cap _{ s'_{\lambda}\in s.'} 
W^{s'_{\lambda}}_{a_{s'_{\lambda}}(J',-1)}  ( C_r^K L) = 
C_r^K (  \cap _{ 
s'_{\lambda}\in s.'} W^{K \cup s'_{\lambda}}_{a_{s'_{\lambda}}(J',-1)}
L)$ is 
defined as above for each subset  of $L$.

\noindent This formula shows that $C_r^{KM}L$ is constructed in two
times, once as $C^K$ over $K_.^+ \times S(K)$ (that is a space normal to $Y_K$) and once 
as a weight filtration over $(M-K)_.^+ \times S(M-K)$ (that is the space $Y_K$).

ii) Let $H' = H^* (C^K_r IC (L )) $ then the above proof apply word for
word to show (notation 18): 
$ C^{KM}_r IC (L ) \simeq {\cal W}_{-1 }\Omega^*(H')$.
For $r \geq 0, H' = 0 $, hence $ C^{KM}_r IC (L )\cong 0 $ and the
isomorphism in (ii) follows since $ C^{KM}_r Q (L ) \simeq  C^{KM}_r  L
/  C^{KM}_r IC (L )$. 

\subsection {Local decomposition.}
\noindent {\em Theorem (decomposition): For a nilpotent otbit $L$ of
dimension 
$n$, there exist canonical injections of $C_r^{KM}L$ in 
$Gr_r^{{\cal W}} (\Omega^{\star} L)$ which decomposes in the category of
perverse  
sheaves (up to a shift in degrees) as a direct sum
$$Gr_r^{\cal W} (\Omega^{\star} L) \simeq \oplus_{K\subset M}C_r^{KM}L.$$
Moreover $ Gr_0^{\cal W} (\Omega^{\star} L) \simeq 0  $ is acyclic.}

\smallskip

\noindent To carry out the proof by induction on $n$, we use only the
property $Gr_0^{\cal W} (\Omega^{\star} L) \simeq 0 $ in dimension $n-1$
to get the 
decomposition in $dim \; n$, then we use the fact that $C_0^K L$ for all
$K$ is acyclic  to get  
again $Gr_0^{\cal W} (\Omega^{\star} L) \simeq 0 $ in dimension $n$ so
to 
complete the induction step. 
For $n = 1$, $K$ and $M$ reduces to one element $1$ and the theorem
reduces to
 $$Gr^W_r \simeq C^1_r L \colon=  Gr^{W^1}_{r+1}L {\buildrel N_1 \over
\rightarrow } Gr^{W^1}_{r-1}L $$ 
By the elementary properties of the weight filtration of $N_1$, it is
quasi-isomorphic to $Gr^{W^1}_{r-1} (L/N_1 L) [-1]$ if $r > 0$,
$Gr^{W^1}_{r+1}(ker N_1:L \rightarrow L)$ if $r < 0 $ and $Gr^{W}_0 L
\simeq  0.$

\noindent For higher dimensions, the proof is carried in various steps.

\smallskip 

\noindent {\em  The complexes $ A^{KM}_r L $, $ B^{KM}_r L$ and
$D^{KM}_r L $}

\smallskip 

\noindent Fixing the $dim.\; n$, the proof is by induction on the length
$\mid K \mid $ of $K$ in $M$. For $r \in \Z$, $K$ fixed and $(J,s.) \in
M_.^+\times S_KM $ we consider the space 
$V(J,s.) =  \bigcap _{s_{\lambda} \varsubsetneq  K, s_{\lambda}\in s.}
W^{s_{\lambda}}_{a_{s_{\lambda}}(J,r-1)} L,$ 
and  the filtrations of $V$: $W^1_t(J,s.) = W^K_{a_K(J,r+ t)} (V(J,s.))
,\;\;
 W^2_t(J,s.) = \bigcap _{K\varsubsetneq s_{\lambda}\in s.} 
W^{s_{\lambda}}_{a_{s_{\lambda}}(J,r+ t)}(V(J,s.)) $
so that 
$W^1_t(J,s.) \cap W^2_t(J,s.)  = \bigcap _{K\subset s_{\lambda}\in s.}
W^{s_{\lambda}}_{a_{s_{\lambda}}(J,r+ t)}L $.

\noindent By summing over $(J,s.)$, we define the complexes

$ A_r^{KM} = s [W^K_{a_K(J,r)}\cap W_0^2 (V (J,s.))/W^K_{a_K(J,r-1)}\cap 
W_{-1}^2 (V (J,s.))]_{(J,s.)\in M.^+\times S_K(M)}  $ 

$B_r^{KM} = s[Gr^{W^2}_0  W^K_{a_K(J,r-1)}( \bigcap _{s_{\lambda} 
\varsubsetneq  K, s_{\lambda}\in s.}
W^{s_{\lambda}}_{a_{s_{\lambda}}(J,r-1)}) 
L]_{(J,s.)\in M.^+\times S_K(M)}  $ 

$ D_r^{KM} = s[Gr^{W^2(J,s.)}_0 Gr^{W^K}_{a_K(J,r)} ( \bigcap
_{s_{\lambda} 
\varsubsetneq  K, s_{\lambda}\in s.}
W^{s_{\lambda}}_{a_{s_{\lambda}}(J,r-1)} 
L)]_{(J,s.)\in
 M.^+\times S_K(M)}  $  

\smallskip 

\noindent { \em Lemma: For all $K \subset M $, there exists a natural
quasi-isomorphism}
$$ B^{KM}_r L  \oplus C^{KM}_r L \cong A^{KM}_r L $$

The proof of this lemma reduces to two sublemmas.

\smallskip

\noindent { \em Sublemma: For all $K \subset M $, there exists an exact
sequence of complexes } $$0\rightarrow B_r^{KM}\oplus
C_r^{KM}\rightarrow  A_r^{KM}\rightarrow  
D_r^{KM}\rightarrow  0$$

The proof is based on the following elementary remark:

\smallskip 

\noindent {\em  Let $W^i$ for $ i =1, 2$  be two  increasing filtrations on an object $V$ 
of an abelian category and $a_i$ two integers, then we have an exact sequence:}

$\quad 0\rightarrow W_{a_2-1}^2 Gr^{W^1}_{a_1}\oplus
W^1_{a_1-1}Gr^{W^2}_{a_2} 
\rightarrow W^1_{a_1}\cap W_{a_2}^2 /W^1_{a_1-1}\cap W_{a_2-1}^2
\rightarrow 
Gr^{W^1}_{a_1}Gr^{W^2}_{a_2}\rightarrow 0 $

\smallskip 

\noindent We apply this remark to the space 
$V(J,s.) =  \bigcap _{s_{\lambda} \varsubsetneq  K, s_{\lambda}\in s.}
W^{s_{\lambda}}_{a_{s_{\lambda}}(J,r-1)} L,$ 
and  the filtrations $W^1$ and $W^2$ of $V$  
for $  a_1 =  0$ and $ a_2 = 0$, then we deduce  from the above sequence
an exact sequence of vector spaces 

$0\rightarrow Gr^{W^2}_0  W^K_{a_K(J,r-1)}(V)\oplus  W^2_{-1}(J,s.)
Gr^{W^K}_{a_K(J,r)}(V) $

$ \rightarrow W^K_{a_K(J,r)}\cap W_0^2 (V )/W^K_{a_K(J,r-1)}\cap
W_{-1}^2 (V) \rightarrow Gr^{W^2}_0 Gr^{W^K}_{a_K(J,r)}(V)\rightarrow 0
$

\noindent The sublemma follows by summing over $(J,s.)$.

\noindent Next we prove by induction on $n$ for the general theorem,
(not only the lemma)

\smallskip 

\noindent { \em Sublemma: $D_r^{KM}\cong 0.$}

\smallskip 

\noindent  Proof. The idea of the proof is to write $D_r^{KM}$ as
$Gr^{\cal W}_{0} (\Omega^{\star} (C^K_r L ))\simeq 0 $ where $ C^K_r L $ is viewed as a 
nilpotent orbit on $M-K$ (that is the fiber of a local system on $Y_{M-K}^*$) and  use 
the induction  
to prove it is zero.
We can either use that $ C^K_r L $ is reduced to its unique non zero
cohomology or as well 
  prove the acyclicity for each term in $C^K_r L $, what we do as follows. 

We simplify the notation from $W^2_t$ above  to

$W^0_t(K,J,s.)\colon=  \bigcap _{K\varsubsetneq s_{\lambda}\in s.} 
W^{s_{\lambda}}_{a_{s_{\lambda}}(J,r+ t)}Gr^{W^K}_{a_K(J,r)} ( \bigcap 
_{s_{\lambda} 
\varsubsetneq  K, s_{\lambda}\in s.}
W^{s_{\lambda}}_{a_{s_{\lambda}}(J,r-1)} L)$ 
\noindent so to write the complex as
$ D^{KM}_r L = s(Gr_0^{W^0 (K,J,s.)} L)_{(J,s.)\in M.^+\times S_K(M)}$, 
then we use the decomposition of $ S_K(M)$ to rewrite $W^0_t$ as

\noindent $ W^0_t(K,(J,J'), (K\cup s.',s.)) =  \bigcap _{s'_{\alpha }\in s.'} W^{K\cup 
s'_{\alpha }}_{a_{s'_{\alpha
}}(J',a_K(J,r)+t)}Gr^{W^K}_{a_K(J,r)} ( 
\bigcap 
_{s_{\lambda} \varsubsetneq  K, s_{\lambda}\in s.} 
W^{s_{\lambda}}_{a_{s_{\lambda}}((J,J'),r-1)} L)$

\noindent $= \bigcap _{s'_{\alpha }\in s.'}  W^{s'_{\alpha
}}_{a_{s'_{\alpha 
}}(J',t)}Gr^{W^K}_{a_K(J,r)} ( \bigcap _{s_{\lambda} \varsubsetneq  K, 
s_{\lambda}\in s.} W^{s_{\lambda}}_{a_{\lambda}(J,r-1)} L)$ 

\smallskip 

\noindent and the complex 

\smallskip 
\noindent $ D^{KM}_r L = s[s[Gr^{W^0}_0(K,(J,J')
(K \cup s.',s.))L]_{(J,s.)\in K.^+\times S(K)}]_{(J',s.')\in
(M-K).^+\times 
S(M-K)}$

\smallskip 
\noindent as a sum in two times over $(J,s.)\in K.^+\times S(K)$ and 
$(J',s.')\in (M-K).^+\times S(M-K))$.
 For a fixed $ (J,s.)$ we consider 
$L(r,J,s.) \colon= Gr^{W^K}_{a_K(J,r)} ( \bigcap _{s_{\lambda}
\varsubsetneq  K, 
s_{\lambda}\in s.} W^{s_{\lambda}}_{a_{s_{\lambda}}(J,r-1)} L)$ 
and the  filtration by subspaces
$ W'^0_t(L(r,J,s.)) \colon= \bigcap _{s'_{\lambda}\in s.'}  
W^{s'_{\lambda}}_{a_{s'_{\lambda}}(J',t)}(L(r,J,s.) $
and finally the complex

\noindent $ D(M-K)(L(r,J,s.) \colon= s[Gr^{W'^0}_0(L(r,J,s.))
]_{(J',s.')\in 
(M-K).^+\times S(M-K)}$.

\noindent By construction
$$ D^{KM}_r L = s[D(M-K)(L(r,J,s.)]_{(J,s.)\in K.^+\times S(K)}.$$
We prove by induction on $n$: $ D(M-K)(L(r,J,s.)\simeq 0. $ 
First we embed $ D(M-K)(L(r,J,s.)) $ in the complex
$$ D'(M-K)(L(r,J)) \colon= s[Gr^{W'^0}_0(Gr^{W^K}_{a_K(J,r)} L
)]_{(J',s.')\in 
(M-K).^+\times S(M-K)} $$
using the embedding  $L(r,J,s.)\subset Gr^{W^K}_{a_K(J,r)} L  $. Now we 
introduce the weight filtration ${\cal W}$ on the combinatorial  DeRham
complex $\Omega^{\star} (Gr^{W^K}_{a_K(J,r)} L)  $ for the nilpotent orbit 
$Gr^{W^K}_{a_K(J,r)} L $ of dimension strictly less then $n$ and weight $a_K(J,r) $
and we notice that  $ D'(M-K)(L(r,J))\simeq Gr^{\cal W}_0
(\Omega^{\star} (Gr^{W^K}_{a_K(J,r)} L ))\simeq 0 $ which is zero by
induction in dimension $n-1$. 
Now $ D'(M-K)(L(r,J))$ is a complex of $MHS$ and  
$W^{s_{\lambda}}_{a_{\lambda}(J,r-1)} $ (up to a shift) is a filtration by subcomplexes 
of $MHS$, so we deduce by strictness that 
for each $r,J,s.$ the complex $ D(M-K)(L(r,J,s.) \simeq 0$ is zero. This ends the proof 
of the sublemma and hence the lemma.

\smallskip 

{ \em Proof of the decomposition theorem}

\smallskip 

\noindent For each $i  \in ${\bf N} we define a map $ \varphi_i: S(M) 
\rightarrow {\cal P}( M)$
to the subsets of $M$ such that  $ M \supset \varphi_i (s.) = Sup
\{s_{\lambda }:
\mid s_{\lambda } \mid \leq i \}$ and for each
 $(J,s.) \in M.^+ \times S(M)$, the filtration with index $t$ of $L$, 
 $W^2_t(\varphi_i (s.),J,s.)\colon=  \bigcap _{\varphi_i (s.)
\varsubsetneq 
s_{\lambda}\in s.} 
(W^{s_{\lambda}}_{a_{s_{\lambda}}(J,r+ t)}) ( \bigcap _{s_{\lambda}
\subset  
\varphi_i (s.) , s_{\lambda}\in s.}
W^{s_{\lambda}}_{a_{s_{\lambda}}(J,r-1)} L)$,
 
\noindent We define  $ G_i(J,s.)(L) = Gr^{W^2(\varphi_i (s.),J,s.)}_0 L$
  then we consider the complex 
$$G_i(s.) = s(G_i(J,s.)L)_{J\subset M}= s(Gr^{W^2(\varphi_i
(s.),J,s.)}_0 
L)_{J\subset M},\; G_i^{SM} L = s(G_i (s.))_{s. \in S(M)}$$
In particular, when $\varphi_i (s.) = \emptyset $, $G_i(s.) =
Gr^{W(s.)}_r L^{s.}$
so that for $i\leq 0$,  $G_i^{SM}L = Gr^{W^{SM}}_r L^{SM}$  and when $i
= \mid M \mid - 1$,  $G_i^{SM}L = C^M_r L$.
Hence the proof of the decomposition theorem will follows from the 

\smallskip 

\noindent {\em Lemma. }
$$G_i^{SM} \cong G_{i+1}^{SM}  \oplus C_r^{KM}$$

Now in order to compare $G_i^{SM} L$ and $G_{i+1}^{SM} L$,
we consider the category  $S^i(M) = \{s. \in S(M): \mid s.\mid \ = i \}$
which is not a 
subcategory of $S(M)$ but $S(M) - S^{i+1}(M)$ is a subcategory such that
the 
restrictions of $G_i$ and $G_{i+1}$ define two subcomplexes:
$G'_i  = s(G_i (s.))_{s. \in S(M)- S^{i+1}(M)}\subset G_i^{SM} L $ 
and $ G'_{i+1} =  s(G_{i+1} (s.))_{s. \in S(M)- S^{i+1}(M)}\subset
G_{i+1}^{SM} L $.
We have $G'_{i+1} = G'_i $ since $\varphi_i  = \varphi_{i+1}$ on $S(M) -
S^{i+1}(M)$.

\noindent The next step is to compute the quotient complexes, for this
we remark:
  $ S^{i+1}(M)\simeq \oplus_{\mid K \mid = i+1}S_K(M)$ and
$\varphi_{i+1} (s.)= K$ for $s. \in S_KM$, then : 
$$G_{i+1}^{S^{i+1}(M)} = G_{i+1}^{SM} L /G'_{i+1} \simeq s(G_{i+1}
(s.))_{s. \in 
S^{i+1}(M)} \simeq  \oplus_{\mid K \mid = i+1} s(G_{i+1} (s.))_{s. \in
S_K(M)}\simeq \oplus_{\mid K \mid = i+1} B_r^KM$$
\noindent On the other side, when $s.\in S_KM$ where $\mid K \mid =
i+1$,
$ G_i(J, s.)=  Gr^{W^1(\varphi_i (s.),J,s.)}_0 L, $ and 
$G_i^{S^{i+1}(M)}  = s(G_i(s.))_{s.\in S_KM}= A_r^KM$, so that
$$ G_i^{S^{i+1}(M)} =   G_i^{SM} L /G'_i \simeq s(G_i (s.))_{s. \in
S^{i+1}(M)} 
\simeq  \oplus_{\mid K \mid = i+1} s(G_i (s.))_{s. \in S_K(M)}\simeq
\oplus_{\mid K \mid = i+1} A_r^KM.$$
Now we deduce from the quasi-isomorphism
$ B_r^KM \oplus C_r^KM \cong  A_r^KM$
that $G_i^{S^{i+1}(M)} = A_r^KM \cong G_{i+1}^{S^{i+1}(M)} \oplus
C_r^{KM}$ hence 
$G_i^{SM}/G_i'\cong (G_{i+1}^{SM}/G_{i+1}')  \oplus C_r^{KM}$, which
proves the lemma since $G_i' \cong  G_{i+1}'  $.

 \subsection { Global construction of the weight filtration.} For each
subset $s_{\lambda} \subset I$ such that $Y_{s_{\lambda} } \neq 
\emptyset$ we write $ {\cal W}^{s_{\lambda}} = {\cal W}(\Sigma_{i\in 
{s_{\lambda}}} {\cal N}_i)$ for the filtration by subbundles defined by
the 
nilpotent endomorphisms of the restriction 
${\cal L}_{Y_{s_{\lambda}}}$ of $ {\cal L}_X$ to $Y_{s_{\lambda}}^*$,
inducing 
also by restriction and for each subset $K \supset  s_{\lambda}$ a
filtration on $ {\cal L}_{Y_K}$. 
Since $({\Omega}^*_X (Log Y) \otimes {\cal L}^{e (\alpha .)}_{X})_y$ is
acyclic if there exists an index $j \in M $ such that $\alpha_j \neq 0 $
(see the formula (11) and the remark below), we can suppose from now on
the local system unipotent.

\smallskip

\noindent {\em Definition ( the weight and  Hodge filtrations). The
weight 
filtration is defined for $ {\cal L} $ unipotent, on the following
combinatorial logarithmic complex $$ {\Omega}^*({\cal L}) =
s({\Omega}^*_{X_{s.}} (Log Y) \otimes 
{\cal L}_{X_{s.}})_{s.\in S}$$
\noindent as follows: let $M \subset I, \mid M\mid = p $ and  $y \in
Y^*_M$, then 
 in terms of a set of $n$ coordinates $y_i, i\in [1,n]$ where we
identify $M$ with $[1,p]$ on an open set $U_y\simeq D^{\mid M\mid}\times
D^{n-p}$ containing $y$ and  a section $f = (f^{s.})_{s.\in S} $ 

$$ f^{s.} = \Sigma_{J\subset M,J'\cap  M = \emptyset} f^{s.}_{J,J'} 
\frac {dy_J}{y_J}\wedge dy_{J'}$$ 
$$ f = (f^{s.})_{s.\in S} \in {\cal W}_r ({\Omega}^*({\cal L}))_{/U_y}
\Leftrightarrow 
\forall J, N \subset M,\, f^{s.}_{J,J'}/Y_N \cap U_y \in {\cap \atop 
{s_{\lambda}\in s.,s_{\lambda} \subset N}} {\cal 
W}_{a_{\lambda}(J,r)}^{s_{\lambda}} ({\cal L}_{/Y_N\cap U_y}) $$

\noindent By convention we let for all integers $r, \;{\cal W}_r / X-Y  = 
{\Omega}^*({\cal L})/ X-Y $, so that ${\cal W}_r$ is the direct image for $r$ big enough 
and the extension by zero for $r$ small enough.
It is a filtration by subcomplexes of analytic subsheaves
globally defined on $X$. The Hodge filtration $F$ is constant in $(s.)$ and deduced from 
Schmid's extension to ${\cal L}_X$
$$F^p (s.) = 0 \rightarrow  F^p{\cal L}_X \ldots \rightarrow  {\Omega}^i_{X_{s.}} (Log Y) 
\otimes 
F^{p-i}{\cal L}_{X_{s.}}\rightarrow \ldots ,\; F^p = s(F^p(s.))_{s.\in
S}$$ }

\noindent { \em Theorem: Consider a unipotent local system $\cal L$ underlying a 
variation of polarised 
Hodge structures of weight $m$; then the complex 
$$ ({\Omega}^*({\cal L}), {\cal W}[m], F)\leqno{(19)}$$ 
with the filtrations ${\cal W}$ and $F$ defined above satisfy 
the decomposition and purity properties. 
More precisely, for all subset $K \subset I$  and all integers $r > 0$ (resp. $ r< 0$),  
let $({\cal L}^K_r, W, F)$ denotes the bifiltered local system underlying a polarised 
$VHS$ on  $Y^*_K $ of weight $r - \vert K \vert + m $ (resp.$r + \vert K \vert + m $) and 
of general fiber $ (H^{\mid K\mid} (C^K_{r} L),W, F)$ of weight induced by  $ W^{s_K}_{r 
- \vert K \vert + m}$ and $F$ defined by $L$
(resp. $ ( H^{\mid K\mid - 1} (C^K_{r} L), W, F)$ for $ r < 0$ of weight induced by  $ 
W^{s_K}_{r + \vert K \vert + m}$ ).
Then we have the following decomposition into intermediate extensions  (up to shift in 
degrees) of (VHS) ${\cal L}^K_r$ compatible with the local decomposition

\smallskip

\centerline {$(Gr^{{\cal W}[m]}_{r+m}{\Omega}^* {\cal L}, F) \simeq \oplus_{K\subset I} 
j_{!*}^K ({\cal L}^K_{r} [-\mid K\mid], W[2\mid K\mid]], F[-\mid K\mid]$, for $ r > 0$ 
and  $j^K:Y^{*}_K \rightarrow Y_K $} 

\smallskip

\centerline {$(Gr^{{\cal W}[m]}_{r+m}{\Omega}^* {\cal L}, F) \simeq \oplus_{K\subset I} 
j_{!*}^K {\cal L}^K_{r} [1 - \mid K\mid], W[-1], F) $, for $ r < 0$}

\smallskip

\centerline {$(Gr^{{\cal W}[m]}_{m}{\Omega}^* {\cal L}, F) \simeq 0$}

\smallskip

that is for $r \geq 0 $ the weight is coincides with the weight for Hodge structures but 
for $r < 0 $ the true  weight for Hodge structures is $ r + m + 1 $ .

iii) The projection on the quotient complex
 $ ({\Omega}^*({\cal L})/j_{!*}{\cal L}, {\cal W}[m], F) $
  with the induced filtrations, induces a filtered quasi- 
isomorphism on $ ( Gr_{r+m}^{{\cal W}[m]}, F)$ for $ r > 0 $.}

\noindent Proof. The decomposition of $(Gr^{{\cal
W}[m]}_{r+m}  {\Omega}^*({\cal L}) ,F) $ reduces near a point $y 
\in Y^*_M $ to the local  decomposition of $ Gr^{\cal W}_{r+m}{\Omega}^* L$ for the 
nilpotent orbit $L$ defined at the point $y$ by the local system
since $C^{KM}_rL $ is precisely the fiber of $j_{!*}^K {\cal L}^K_{r} [-\mid K\mid]$ for 
$r >0$ (resp. $j_{!*}^K {\cal L}^K_{r} [1 -\mid K\mid]$
for $r < 0$. The count of weight takes into account for $r > 0$ the residue in the 
isomorphism with $L$ that shifts $W$ and $ F$ but also the shift in degrees, while for $ 
r < 0 $ there is no residue  but only a shift in degrees, the rule being as follows:

{\em Let $( K, W, F)$ be a mixed Hodge complex then for all $m, h \in \Z$, $( K', W', F') 
= ( K [m], W[m-2h], F[h])$ is also a mixed Hodge complex.}

   The same proof apply for $r =0$, hence 
 ${\cal W}_{- 1} \simeq {\cal W}_0 $ is isomorphic to the intermediate 
extension of ${\cal L}$ by Kashiwara and Kawai's formula, that we prove below. The 
assertion (iii) follows from the
assertion (iii)  in the purity theorem corresponding to a result on $C^{KM}_r QL $.

\smallskip

 {\it  Proof of Kashiwara and Kawai's formula: $j_{!*} {\cal L}[2n]
\simeq {\cal W}_0{\Omega}^*({\cal L}[2n])$.}

\smallskip
 
\noindent In this subsection we give a proof of the formula of the
intermediate extension of ${\cal L}[2n]$,  announced in  [26], which is
in fact the subcomplex ${\cal W}_0{\Omega}^*({\cal L}[2n])$. It follows
easily from the local decomposition of the graded parts of the weight
filtration, by induction on the dimension $n$.

\smallskip
\noindent {\em Theorem. The subcomplex $ {\cal W}_0{\Omega}^*({\cal
L}[2n])$ is quasi-isomorphic to the intermediate extension of ${\cal
L}[2n]$.}

\smallskip
The proof of this theorem is by induction on  the dimension $n$. It is
true in dimension $1$ 
and if we suppose the result true in dimension strictly less than $n$, 
we can apply the result for local systems defined on open subsets of the
closed sets $Y_K$, namely the local system ${\cal L}^K_r [-\mid K\mid])$
for $r > 0$ (resp.${\cal L}^K_r [1 -\mid K\mid])$
for $r < 0$) 
whose fiber at at  each point $y \in Y^*_K$ is quasi-isomorphic to
$C^K_rL$. Let $ j^K:Y^{*}_K \rightarrow Y_K $ be  the open embedding in
$Y_K$ and consider the associated DeRham complex 
${\Omega}^*({\cal L}^K_r)$ 
on $Y_K$ whose  weight filtration will be denoted locally near a point in 
$Y^*_M$ by ${\cal W}^{M-K}$ for $K \subset M$; then by the induction hypothesis we have 
at 
 the point $y $:  ${\cal W}_{-1} {\Omega}^* L \simeq {\cal W}_0
{\Omega}^* L  $ is also quasi-isomorphic to the fiber of the
intermediate extension of ${\cal L}$, that is
 $$ \forall r > 0, C_r^{KM}(L)\simeq (j_{!*}^K {\cal L}^K_{r} [-\mid K\mid])_y \simeq 
{\cal W}^{M-K}_{-1 } C^K_r L$$ 
and similarly for $r < 0$.

We use the following criteria caracterising  intermediate extension
[17]:

\noindent Consider the stratification defined by $Y$ on $X$ and the
middle perversity
$p(2k) = k -1 $ associated to the closed subset $Y^{2k} = \cup_{\mid
K\mid = k} Y_K$ of real codimension $2k$. We let $Y^{2k-1} = Y^{2k} $
and $p(2k-1) = k -1 $. For any complex of sheaves $S.$ on $X$ which is
constructible with respect to 
the stratification , let $S.^{2k} = S.^{2k-1}  = S.\mid X-Y^{2k}$ and
consider the four properties:

a) Normalisation: $S.\mid X-Y^2 \cong {\cal L}[2n]$

b) Lower bound: ${\bf H}^i (S.) = 0 $ for all $i < -2n$

c) vanishing condition: ${\bf H}^m ( S.^{2(k+1)})= {\bf H}^m ( S.^{2k
+1})= 0  $ for all $m > k -2n$

d) dual condition: 
${\bf H}^m (j_{2k}^! S.^{2(k+1)})= 0  $ for all $k \geq 1$ and all $m >
k -2n$ where  $j_{2k} : Y^{2k} - Y^{2(k+1)}\rightarrow X - Y^{2(k+1)}$
is  the closed embedding,

 \noindent then $S.$ is the intermediate extension of ${\cal L}[2n]$.
 
\noindent In order to prove the result for $n$ we check the above four
properties for $ {\cal W}_0{\Omega}^*({\cal L}[2n])$. The first two are
clear and  we use the exact sequences $$ 0 \rightarrow
{\cal W}_{r-1}\rightarrow {\cal W}_r \rightarrow Gr^{\cal W}_r
\rightarrow 0 $$ 
to prove d)(resp. c)) by descending (resp. ascending )indices from
${\cal W}_r $ to ${\cal W}_{r-1}$ for $r \geq 0$ (resp. $r-1 $ to $r$ for $r < 0 $ ) 
applying at each step the inductive hypothesis to $Gr^{\cal W}_r $.

Proof of d). The dual condition is true for $r$ big enough since then $
{\cal W}_r $ coincides with the whole complex, that is the higher direct
image of ${\cal L}[2n]$ on $X-Y$. Now we apply d) on $Y_{K'} $ with
$\mid K'\mid \,= \,k'$  for  $j'_{2k} : Y^{2k}\cap Y_{K'}  - Y^{2(k+1)}\cap
Y_{K'} \rightarrow Y_{K'}  - Y^{2(k+1)}\cap Y_{K'} $ where we suppose $k
> k'$ ( notice that $Y^{2k}\cap Y_{K'} = (Y \cap Y_{K'})^{2(k-k')}$, then for $S.' $ 
equal to the intermediate extension of ${\cal L}^{K'}_r [2n - 2k']$ on $Y_{K'}$ we have  
the property ${\bf H}^m ({j'}^!_{2k} S.'^{2(k-k')+1)})= 0
$ for all $ (k - k') \geq 1$ and all $m > k - k' - 2(n-k') = k+k' - 2n$ which gives for $
S.' [k'] $ on X: ${\bf H}^m (j_{2k}^! S.'^{2(k+1)}[k'] )= 0  $ for all $k > k' $ and all 
$m > k  - 2n$, hence d) is true.

 If $k = k'$, then $Y^{2k}\cap Y_{K'}  =
Y_{K'}$ and we have a local system in degree $k'-2n $ on $Y_{K'}  - Y^{2(k+1)}\cap Y_{K'} 
$ hence d) is still true
and for $k < k' $, the support   $Y_{K'}  - Y^{2(k+1)}\cap Y_{K'} $ of $S.'$ is empty. 
From the decomposition theorem and the induction, this
argument apply to $Gr^{\cal W}_r $ and hence apply by induction on $r \geq 0$ to $ {\cal
W}_0 $ and also to ${\cal W}_{-1} $.

Proof of c). Dually, the vanishing condition is true for $r$ small
enough since then $  {\cal W}_r $ coincides with the  extension by zero
of ${\cal L}[2n]$ on $X-Y$.

\noindent Now we use the filtration for $r < 0 $, for $S.' $ equal to the intermediate
extension of ${\cal L}^{K'}_r [2n - 2k']$ on $ Y^{K'}$ we have for $k > k'$: ${\bf H}^m 
(S.'^{2(k-k')+1})= 0  $ for all 
$m > k + k' - 2n$, which gives for $ S.' [k'+1],\; (r < 0) $ on X: ${\bf H}^m ( 
S.^{2(k+1)})= {\bf H}^m ( S.^{2k +1})= 0  $ for all $m > k-1 -2n$. If $k = k'$, then $ 
S.' [k'+1] $ is a local system in degree  $-2n + k-1 $ on $Y_{K'} - Y^{k+1}$ and for $k < 
k' $, $Y_{K'} - Y^{k+1}$ is empty.

\smallskip
\noindent { \em Corollary: If we suppose $X$ proper and we replace the filtration ${\cal 
W}$
by ${\cal W}''$ with ${\cal W}_i'' = {\cal W}_i $ for $i \geq 0$
and ${\cal W}_{-1}'' = 0$,  then the bifiltered complex 
$$ ({\Omega}^*({\cal L}), {\cal W}''[m], F)$$ 
is a  mixed Hodge complex .} 

\section{The complex of nearby cycles  ${\Psi}_f (\cal L)$.} 
Let $f \colon X \rightarrow D$ and suppose $Y = f^{-1}(0)$; the definition of 
the complex of sheaves of nearby cocycles on $Y$ is given in [11]; its 
cohomology fiber at a point $y$ equals the cohomology of the Milnor fiber $F_y$  at $y$ 
in $Y$. 
The monodromy ${\cal T}$ induces an action on the cohomology 
$H^i ({\Psi}_f {(\cal L)}_y) \simeq H^i (F_y , {\cal L})$ and on the complex itself 
viewed in the abelian 
category of perverse sheaves. It is important to point out that the action on  the 
complex is related to the action on cohomology through a spectral sequence 
and precisely in our subject we need to use the weight  filtration on the  complex itself 
and not on its cohomology.
 
\noindent  The aim of this section is to describe the weight filtration
on ${\Psi}_f ({\cal L})$. This problem is closely related to the weight
filtration in the open case since there exists a close {\it relation
between  ${\Psi}^u_f ({\cal L})$,  the direct image $ {\bf j}_* {\cal
L}$ and ${\bf j}_{!*} {\cal L}$ } as  explained in [2] ( and previously in a private 
letter by Deligne and Gabber)

\smallskip
\noindent {\em  Proposition [2]: Let ${\cal N} = Log{\cal T}^u$ denotes the logarithm of 
the unipotent
part of the monodromy, then we have the following isomorphism in the abelian category of 
perverse sheaves
$${\bf j}_* {\cal L} / {\bf j}_{!*} {\cal L} \quad \simeq \quad Coker
({\cal N} 
\colon {\Psi}_f  ^u {(\cal L)} \rightarrow 
{\Psi}_f^u {(\cal L)}) \lbrack -1 \rbrack \leqno{(20)} $$ }

\noindent The  filtration $W({\cal N})$ on ${\Psi}_f^u
({\cal L})$  induces a filtration ${\cal W}$ on $Coker {\cal N} 
/ {\Psi}_f^u ({\cal L})$, hence on ${\bf j}_* {\cal L} / {\bf j}_{!*}
{\cal L}$.

\smallskip

{\it  The induced  filtration on 
${\bf j}_* {\cal L} / {\bf j}_{!*}  {\cal L}$  is  independant of the choice of $f$.} 
For a rigorous proof one should use the result of Verdier [34].
To prove  the independance of $f$ we can use a path in the space of 
functions between two local equations $f$ and $f'$ of $Y$ and defines
by parallel  
transport an  isomorphism between  ${\Psi}_f ({\cal L})$ and 
${\Psi}_{f'} ({\cal L})$ ; modulo $coker {\cal N}$, this isomorphism is
independant of the path.

\subsection{ The weight filtration on the nearby cycles ${\Psi}_f ({\cal
L})$}

The method to compute ${\Psi}_f$ as explained in [11] uses the
restriction 
$i^*_Y \hbox {\bf j}_* {\cal L}$ of the higher direct image  of ${\cal
L}$ to $Y$ and the cup-product $H^i(X^*, {\cal L}) \otimes H^1(X^*, \Q ) \;
{\buildrel {\smile  \eta } \over \longrightarrow} \;H^{i+1}(X^*, {\cal
L})$ by the inverse image 
$\eta = f^* c \in H^1 (X^*, \Q )$ of a  generator $c$ of the 
 cohomology $H^1(D^*, \Q )$.
Thus one defines a morphism (of  degree 1),
$\eta \colon \; i^*_Y \hbox {\bf j}_* {\cal L} \to i^*_Y 
\hbox {\bf j}_* {\cal L} [1 ]$
such that  ${\eta}^2 = 0$ so to get a double complex whose simple
associated 
complex is quasi-isomorphic to  ${\Psi}^u_f ({\cal L})$, the unipotent part of  
${\Psi}_f ({\cal L})$ under the monodromy action ${\cal T}$
$${\Psi}^u_f ({\cal L})\;
\simeq \; s{(i^*_Y \hbox {\bf j}_* {\cal L} [ p ], \eta)}_{p \leq 0} \leqno{(21)}$$
In order to get the full ${\Psi}_f ({\cal L})$ (not only the unipotent part under the 
action of ${\cal T}$) Deligne introduced local systems of rank
one $ {\cal V}_{\beta }$ on the disc with monodromy $e(\beta)= exp (-2i\pi \beta 
)$ and proved the following isomorphism 
$$ \Psi_f ({\cal L})\; \simeq \; \oplus_{\beta \in {\C}} \; \Psi^u_f({\cal L} \otimes 
f^{-1}{\cal V}_{\beta })\leqno{(22)} $$
When   ${\cal L}$ is quasi-unipotent we need only to consider $\beta \in {\Q} \cap 
[0,1[$. Moreover, near a point $ y \in Y$ such that $f = \Pi_{j\in M}  z_j^{n_j}$, the 
tensor product ${\cal L}^{e(\alpha .)} \otimes f^{-1}{\cal 
V}_{\beta } $ is unipotent near $y$ if and only if 
$ \forall j \in M, \alpha_j + n_j \beta \in {\bf N}$, then
 $$ \Psi_f ({\cal L}^{e(\alpha .)}) \; \simeq \; \oplus_{\beta \in S} \; \Psi^u_f ({\cal 
L}^{e(\alpha .)}\otimes f^{-1}{\cal V}_{\beta }),\quad S = 
\{\beta   \in {\C}: \forall j \in M, \alpha_j + n_j \beta  \in \hbox {\bf 
N}\}.\leqno{(23)}  $$
Due to this formula, the problem can be reduced later in  the article to study
the unipotent part  $\Psi^u_f$. We recall that in order to construct ${\cal L}_X$ we need 
to choose a section as follows

\smallskip

 \noindent {\em Definition: We define $\tau:\C/\Z \to \C $ as the section of 
$\pi : \C \to \C/\Z $ such that $ Re (\tau) \in [0,1[. $ }

\smallskip

\noindent {\it  Local description }.
Near a point $y \in Y$, where $f = {n \atop {\buildrel \prod \over {i=1}}} 
z^{n_i}_i$ for non zero $n_i$, in  DeRham cohomology  $\eta =
f^*(\frac{dt}{t}) = {n \atop {\buildrel \Sigma \over {i=1}} }  n_i \frac {dz_i} {z_i}$. 
The morphism $\eta$ on 
$\Omega (L,D. + N.)$
(3) is defined by $ {\varepsilon}_{i_j} n_{i_j} Id: L(i. - i_j) = L \rightarrow  L = 
L(i.)$
where ${\varepsilon}_{ij}$ is the signature of the permutation which order strictly $( i. 
- i_j, i_j )$ 
for various $i_j$. For each complex
number $\beta$, we consider the following complexes where 
$L^{e(\alpha.)}\; = \; \cap_{i\in [1,n]} L^{e(\alpha_i)} $ is the
intersection of the eigenspaces for  $T_i^s, i\in [1,n]$  with
eigenvalues $ e(\alpha_i)$. $${\Psi}^{\beta}_p (L)\; = \; \oplus_{\alpha
.} \; \Omega (L^{e(\alpha.)} , \tau ( \alpha_i + \beta n_i ) Id + N_i
)_{i\in [1,n]})\; [p ], \;p \leq 0 
\leqno{(24)} $$
where $\eta \colon {\Psi}^{\beta}_p (L) \rightarrow
{\Psi}^{\beta}_{p+1}(L)$ is a complex morphism
satisfying  ${\eta}^2 =0$, the $({\Psi}^{\beta}_p,\eta ) $ form a double complex 
for $p \leq 0$. Let ${\Psi}^{\beta}(L)$ denotes the  associated simple complex. In order 
to take into account the action of  
$ N = - 1/2i\pi Log T^u$ we  write after Kashiwara,
$L[N^p]$ for $ L[p]$ and $ L [N^{-1}]$ for the direct sum over $p$, so that the 
action of $ N $ is just multiplication by $N$ 
 $${\Psi}^{\beta}(L) = {s({\Psi}^{\beta}_p(L),\eta )}_{p \leq 0} \simeq \oplus_{\alpha .} 
\Omega (L^{e(\alpha.)}[N^{-1} ],  \tau ( \alpha_i + \beta n_i ) Id + N_i - n_i N)_{i\in 
[1,n]}). \leqno{(25)} $$

It is isomorphic to the direct sum of Koszul complexes defined by $
(L^{e(\alpha.)}[N^{-1} ],  \tau ( \alpha_i + \beta n_i ) Id + N_i - n_i N)_{i\in [1,n]}. 
$ 
The complex $\Omega (L^{e(\alpha .)} , \tau (\beta n.+ \alpha .) Id + N.)$ is acyclic 
unless
 $ (n_i\beta + {\alpha}_i) \in {\bf N} $ for all $i \in M$, hence 
${\Psi}^{\beta}(L)$ is acyclic but for a finite number of $\beta$ such that $e(\beta)$ is 
an eigenvalue  of the  monodromy action. The proof of Deligne's result [11] reduces to:
 
\smallskip 
\noindent {\it The fiber  at zero in $D^{n+k}$ of ${\Psi}_f (\cal L)$ (resp. 
${\Psi}_f^u (\cal L))$ is quasi-isomorphic to  a (finite )  direct sum of $\Psi^{\beta} 
(L)  (25) $(resp. to  ${\Psi}^0(L) $ for $\beta = 0$)}

$${\Psi}_f {(\cal L)}_0 \; \simeq \;\oplus_{ \beta \in \C }
{\Psi}^{\beta}  L  
\quad , \quad {\Psi}^u_f {(\cal L)}_0 
\simeq {\Psi}^0 L  \leqno{(26)}$$

\be{The weight and Hodge filtrations on ${\Psi}^{0}L $}{3.I.1} \
We consider again a nilpotent orbit $L$.
To describe the weight in terms of the filtrations $ ({\Omega}^* L, {\cal W}, F)$ 
associated to $L$, we need to use the constant complex with index $s. \in S(M)$, 
${\Psi}^{0} L(s.) = {\Psi}^{0} L $ and  introduce the complex 
$$ ({\Psi}^{0} L)_{ S(M)} \colon= s({\Psi}^{0} L (s.))_{s. \in S(M)} \leqno{(27)} $$
 which can be viewed also as $ s({\Omega}^* L [p], \eta)_{p \leq 0}$, then we define on 
it the weight filtration
$$ {\cal W}_r ({\Psi}^{0} L)_{ S(M)} =  s ({\cal W}_{r+2p - 1} {\Omega}^* L [p], \eta)_{p 
\leq 0}, \; \; 
F^r ({\Psi}^{0} L)_{ S(M)} =  s ( F^{r+p} {\Omega}^* L [p], \eta)_{p \leq 0}.  
\leqno{(28)}$$
\noindent {\em Monodromy.} 

 The logarithm ${\cal N}$ of the monodromy is defined by an endomorphism
$\nu$ of the  complex $ {\Psi}^0(L)_{ S(M)} $, given by the formula

$$\forall a. = \Sigma_{p\leq 0} a_p \in  ({\Psi}^0 L)_{ S(M)},(\nu
(a.))_p = a_{p-1} \; \forall p \leq 0$$
such that $\nu ( {\cal W}_r ) \subset {\cal W}_{r-2}$ and  
$\nu ( F^r ) \subset F^{r-1}$.

\smallskip 

\noindent {\em Decomposition of} $ Gr^{\cal W}_r$ 

The morphism $\eta $ induces a morphism denoted also by $\eta : C^{KM}_r
L \rightarrow C^{KM}_{r+2} L [1]$ so that we can define a double complex
and the associated simple complex

$$ \Psi^{KM}_r L = s( C^{KM}_{r+2p-1} L [p], \eta )_{p \leq 0}, \;
\Psi^{K}_r L \colon= \Psi^{KK}_r L $$  
We will see soon that this complex decomposes into a direct sum.

\smallskip 
\noindent {\em Lemma:
There exists  natural injections of $\Psi^{KM}_r L $ into 
$ Gr^{\cal W}_r ({\Psi}^0 L)_{ S(M)} $ and a decomposition} 
$$  Gr^{\cal W}_r ({\Psi}^{0} L)_{ S(M)}  \simeq \oplus_{K \subset M}
\Psi^{KM}_r L \leqno{(29)}$$
Proof: By the spectral sequence of a double complex, it is enough to
check the decomposition on the columns where the proof reduces to the
decomposition in the open case.

\smallskip 
\noindent {\em Theorem: 
\noindent The weight filtration  (28) coincides with $W({\cal N})$ defined by the 
logarithm of the monodromy in the abelian category of perverse sheaves.}

\smallskip 
\noindent The proof in two steps reduces to the lamma and 
the proposition below.

\smallskip 

\noindent {\em Lemma:
The following statements are equivalent

\noindent i) For all $i\geq 1, \; \nu^i :   Gr^{\cal W}_i ({\Psi}^0 L)_{
S(M)} \simeq  Gr^{\cal W}_{-i} ({\Psi}^0 L)_{ S(M)}$.

\noindent ii)  For all $i\geq 1 , \; Gr^{\cal W}_i ker \; \nu^i =
 s(Gr^{\cal W}_{i+2p-1}{\Omega}^* L [p], \eta)_{-i <p \leq 0} \cong 0$}.

\smallskip 

\noindent Proof: 
 The morphism $\nu^i$ on 
$({\Psi}^0 L)_{ S(M)} $ is surjective and its kernel is sum of the
columns of index $ -i < p \leq 0 $.

\smallskip 

{\em Remark}: It may be interesting for the reader to check the statement on the example 
of a line with $f$ equivalent at $0$ to $z^n$ on the fiber of $({\Psi}^0 L)$ at the point 
$0$ for $ L = \C $ and $N=0$, where the similarity and the differences with Steenbrink's  
construction  appears already.

\smallskip
\noindent {\em Proposition:
 For all $i \geq 1, \; Gr^{\cal W}_i ker \; \nu^i =  s(Gr^{\cal
W}_{i+2p-1}{\Omega}^* L [p], \eta)_{-i <p \leq 0} \cong 0$.}

\smallskip 
 Proof: Let $ \nu^{KM}_i \colon= s[C^{KM}_{i+2p-1} L [p], \eta]_{-i < p
\leq 0}$, then by the decomposition theorem we have: $ Gr^{\cal W}_i ker
\; \nu^i \cong \oplus_{K \subset M} \nu^{KM}_i .$ 
Denotes $ \nu^{KK}_i$ by $ \nu^{K}_i$, then we can easily check that $
\nu^{KM}_i$ is the intermediate extension of $ \nu^{K}_i$ and is
quasi-isomorphic to zero if $ \nu^{K}_i \cong 0$, so we reduce the proof
to

\smallskip 

\noindent {\em Lemma: For all $i \geq 1,\; \nu^{K}_i \colon=
s[C^K_{i+2p-1} L [p], \eta]_{- i < p \leq 0} \cong 0$}

Proof: In order to give a proof  by induction for $i$ assuming the
result for $i-2$ , we write $ \nu^{K}_i$ as:

$s [ C^K_{-(i-1)} L [-(i-1)], s[C^K_{i+2p-1} L [p], \eta]_{-(i-2) < p
\leq 0} [-1], C^K_{i-1} L, \eta] \cong 0 $

\noindent We know that 
$  C^K_{-(i-1)} L [\mid K \mid -1] \cong Gr^{W^K}_{\mid K \mid
-(i-1)}[(\cap_{i \in K} ker N_i: 
L \rightarrow L ] \simeq $

$\oplus_{(m_1,\cdots,m_n) \in T'(r)} Gr^{W^n}_{m_n } \cdots 
Gr^{W^i}_{m_i } \cdots Gr^{W^1}_{m_1 } [(\cap_{i \in K} (ker N_i: 
L \rightarrow L ] $

\noindent and
$ C^K_{(i-1)} L [\mid K \mid ] \cong Gr^{W^K}_{r- \mid K
\mid}[L/(\Sigma_{i \in K} N_i L)]$

$ \simeq \oplus_{(m_1,\cdots,m_n) \in T(r)} Gr^{W^n}_{m_n - 2 } \cdots 
Gr^{W^i}_{m_i - 2 } \cdots Gr^{W^1}_{m_1 - 2 } [L/(\Sigma_{i \in K} N_i
L)] $.

\noindent Given a nilpotent orbit $(L, N_i)$, we denote in general the
primitive part of $(Gr^{W(N)}_{r} L)$ by $(Gr^{W(N)}_{r} L)^0$, then we
have the following isomorphisms :

 $Gr^{W(N)}_{r} (L/N L) \simeq (Gr^{W(N)}_{r} L)^0
{\buildrel { N^r} \over {\rightarrow }} (Gr^{W(N)}_{- r} ker N)$,

\noindent  so we can deduce in general:

 $ N_1^{m_1} \cdots N_n^{m_n}: Gr^{W^n}_{m_n } \cdots 
Gr^{W^i}_{m_i } \cdots Gr^{W^1}_{m_1 } [L/
(\Sigma_{i \in K} N_i L)] \simeq  Gr^{W^n}_{m_n } \cdots 
Gr^{W^i}_{m_i } \cdots Gr^{W^1}_{m_1 } \cap_{i \in K} ker N_i$

\noindent the sum over
 $\{ ( m_1 \geq 0, \cdots, m_n \geq 0 ) \}:
 (\Sigma_{i \in K} m_i = i-1 - \mid K \mid $ induces an isomorphism:

$ \gamma: Gr^{W^K}_{i-1 - \mid K \mid}[L/
(\Sigma_{i \in K} N_i L)] \rightarrow  Gr^{W^K}_{-(i-1)+\mid K
\mid}[(\cap_{i \in K} ( ker N_i: L \rightarrow L ] $

\noindent  since 
$ Gr^{W^K}_{-(i-1)+\mid K \mid}[\cap_{i \in K} ker N_i: 
L \rightarrow L ] $  is isomorphic to

$ \oplus_{\{ ( m_1 \leq 0, \cdots, m_n \leq 0 ) \}:
 (\sum_{i \in K} m_i =  \mid K \mid - i + 1} \cdots \oplus Gr^{W^n}_{m_n
} \cdots Gr^{W^i}_{m_i } \cdots Gr^{W^1}_{m_1 }$.
 
 \noindent Then $ \gamma$ induces a quasi-isomorphism from

$ C^K_{(i-1)} L [\mid K \mid ] \cong Gr^{W^K}_{i-1- \mid K
\mid}[L/(\Sigma_{i \in K} N_i L)] \cong Gr^{W^K}_{i-1 - \mid K
\mid}[(L/( N_1 L)/N_2 (L/( N_1 L)] $ 

\noindent to $ C^K_{(1 -i)} L [\mid K \mid -1 ] $.
 
\noindent A diagram chasing shows that  $\nu^{K}_i$ is in fact a cone
over $ \gamma ^{-1}$, hence zero, which establishes the lemma and the
proposition.

\smallskip 
\noindent {\em Corollary: The  graded part of $({\Psi}^0 L)_{ S(M)}$
is non zero for only a finite number of indices for which it reduces to a 
double 
complex of finite terms. More precisely, let  $i_0 $ be an integer large enough 
to have $ Gr^{\cal W}_{j}{\Omega}^* {\cal L} = 0$ for all $ j > i_0$,then for a given $i 
\geq 1$:

$$Gr^{\cal W}_{i}({\Psi}^0 L)_{ S(M)} \cong  s(Gr^{\cal W}_{i+2p-1}{\Omega}^* L 
[p], \eta)_{ p \leq - i}$$
 where only a finite number of $p$ such that
$ - i_0 \leq i + 2p -1 \leq -i-1$ are non zero in the right term. For 
$i \leq - 1$ we use the isomorphism $\nu^{-i}$.
In particular, $Gr^{\cal W}_{i}({\Psi}^0 L)_{ S(M)} \cong 0$
for all $i$ such that  $ \vert i \vert \geq i_0 $. We have a direct sum of intermediate 
extensions (up to shift in degrees) of $VHS$ of  weight $i + m + 1 $ }. 

\smallskip 

\noindent Proof: Suppose $i > 0$, then $Gr^{\cal W}_{i}({\Psi}^0 L)_{ S(M)} =
  s(Gr^{\cal W}_{i+2p-1}{\Omega}^* L [p], \eta)_{ p \leq 0} $ is the cone
over 

 $ \eta: s(Gr^{\cal W}_{i+2p-1}{\Omega}^* L [p], \eta)_{ -i < p
\leq 0}[-1] \rightarrow   s(Gr^{\cal W}_{i+2p-1}{\Omega}^* L [p], \eta)_{ p \leq 
-i} $

\noindent where the first complex is 
$ Gr^{\cal W}_i ker \; \nu^i $ hence quasi-isomorphic to zero, then the corollary 
follows.

\smallskip 

\noindent {\em Remark: i) This corollary, shows that the weight
filtration behaves like a finite one, so that we can apply in the proper
case the results on mixed Hodge complex where the weight filtration is
 supposed to be finite.}

\smallskip 

\noindent {\em ii)
In Steenbrink's case that is ${\cal L} = \C$, 
 $ s({\cal W}_{p-1}{\Omega}^* L [p], \eta)_{-{i_0} < p \leq
0}$ is a subcomplex quasi-isomorphic to
$({\Psi}^0 L)_{ S(M)}$. For a general ${\cal L}$, it is not a complex, 
nevertheless the graded part behaves like if we restrict to such object.}

\smallskip 

\noindent {\em iii) Dually, we could define 
$ ({\Psi}^{0} L)_{ S(M)} $ as $i_Y^{\s} s({\Omega}^* L [p], \eta)_{p \geq 0}[1]$, with 
the filtrations
$$ {\cal W}_r =  s ({\cal W}_{r+2p+1} {\Omega}^* L [p], \eta)_{p \geq 0}[1],
\; F^r =  s ( F^{r+p+1} {\Omega}^* L [p], \eta)_{p \geq 0}[1]$$
then the above results show that the two definitions give quasi-isomorphic 
complexes and the formula for $ Gr^{\cal W}$ behaves like if we could use the quotient by 
${\cal W}_p\Omega^* L $ in each column $p$.

We will see later that we can take the quotient with the subcomplex generated by $IC L $ 
for $p = 0$ and then use the induced filtrations on the quotient.}

\smallskip 

\noindent {\em Corollary (decomposition): Let $ I(p) = 
\{ p \leq 0 ,  - i_0 \leq i + 2p -1 \leq -i-1 \}$, then:}
$$Gr^{\cal W}_{i}({\Psi}^0 L)_{ S(M)} = 
 s(Gr^{\cal W}_{i+2p-1}{\Omega}^* L [p], \eta)_{ p\in I(p)} = 
 \oplus_{p \in I(p)} Gr^{\cal W}_{i+2p-1} {\Omega}^* L [p].$$
\noindent Proof : It follows from the remark that $Gr^{\cal W}_{i}({\Psi}^0 L)_{ S(M)} $ 
can be computed for a finite number of columns such that $r = i+ 2p - 1 \leq p-1 < 0 $
is negative and where each term is a direct sum of $C^{KM}_r L$, intermediate extension 
of $C^{K}_r L$ whose  cohomology is concentrated in degree 
$\vert K \vert - 1$, hence the map $\eta $ is zero and we get a direct sum instead of a 
double complex.

\be{The global weighted complex  $({\Psi}^u_f ({\cal L}), {\cal W}, F)$}
{3.I.2} \ 

Returning to the global situation, we need to define the Hodge filtration on 
$\Psi_f ({\cal L}_X)$. First $F$ extends to the logarithmic complex by the formula: 
$ F^p ({\Omega}^*_X (Log Y) \otimes {\cal L}_X) = 
s ( {\Omega}^q_X (Log Y) \otimes F^{p-q} ({\cal L}_X), \nabla_X)_{p \leq 0}) $, 
then $F$ extends to $(\Psi_f ({\cal L}_X)$  via the formula 
 $$ F^r( i_Y^{\s} s({\Omega}^*_X (Log Y) \otimes {\cal L}_X[i])_{i\leq 0}) =     i_Y^{\s} 
s(F^{r+p+1}(\Omega^*_X (Log Y) \otimes {\cal L}_X ) [i], \eta)_{i\leq 0}$$
 
The definition of the global weight filtration reduces to the local construction at a 
point  $y \in Y^*_M$, using the quasi-isomorphism
${(\Psi}_f \;({\cal L}^{e(\alpha.)}))_y \; \simeq \;\oplus_{ \beta} \;{\Psi}^{\beta 
}(L^{e(\alpha .)})$.  
	
We suppose again ${\cal L}$ unipotent and define 
as previously 
$ ({\Psi}^{u}_f {\cal L})_{ S(M)} \colon= s({\Psi}^{u}_f {\cal L}(s.))_{s. \in S(M)}$
 which can be viewed also as $s({\Omega}^* {\cal L} [p],
 \eta)_{p \leq 0}$,
then we define on it the weight filtration
$$ {\cal W}_r ({\Psi}^{u}_f {\cal L})_{ S(M)} = i_Y^{\s} s ({\cal W}_{r+2p-1} {\Omega}^* 
{\cal L} [p], \eta)_{p \leq 0}, \; 
F^r ({\Psi}^{u}_f {\cal L})_{ S(M)} = i_Y^{\s} s ( F^{r+p+1} {\Omega}^* {\cal L}[p], 
\eta)_{p \leq 0}\leqno{(30)}$$
The logarithm of the monodromy ${\cal N}$ is defined on this complex as
in the local case.
The filtration $W({\cal N})$ is defined on  $({\Psi}^{u}_f {\cal L})_{
S(M)}$ in the abelian category 
of  perverse sheaves.

\smallskip

\noindent { \em Theorem: Suppose $\cal L$ underlies a unipotent variation of polarised 
Hodge structures of weight $m$, then the graded part of the weight filtration (30) of the 
complex }

\centerline {$(\Psi^u_f ({\cal L}_X),{\cal W}[m], F)$}
 
\noi { \em decomposes into a direct sum of intermediate extension of VHS; moreover we 
have $W({\cal N})= {\cal W}$.}

\smallskip

The proof of this theorem reduces by definition to show that 
$(Gr^{{\cal W}[m]}_r,F) $ decomposes which result can be reduced to 
the local case where it has been checked in the above corollaries. 

\smallskip

\noindent { \em Remark: We could as well  define the complex 
$(\Psi^0 L)_{S(M)}$ by summing over $ p \geq 0 $:
 $$ (\Psi^0 {\cal L })_{S(M)} = i_Y^{\s} s[(s(\Omega^*_X (Log Y) \otimes 
 {\cal L}_X )[p], \eta)]_{p\geq 0}[1] \leqno{(31)}$$
 
 By the above remarks the two definitions give quasi-isomorphic complexes.}

\be{The global weighted complex of $({\Psi}_f (\cal L), W, F)$}{3.I.3} \
Let $y \in Y^*_M$, we deduce from the
isomorphism  $${(\Psi}_f \;({\cal L}^{e(\alpha.)}))_y \; \simeq
\;\oplus_{ \beta} \;{\Psi}^{\beta }(L^{e(\alpha .)})\leqno{(32)} $$  

\noindent the global weight filtration in the abelian category 
of  perverse sheaves on  Deligne's extension  $({\cal L}^{e (\alpha .)}\otimes 
f^{-1}{\cal V}_{\beta })_X$ and  the associated combinatorial 
 logarithmic complex
$ {\Omega}^* ({\cal L}^{e (\alpha.)}\otimes f^{-1}{\cal V}_{\beta })$
where we define the global weight filtration in the abelian category 
of  perverse sheaves at points $y$ such that ${\cal L}^{e (\alpha.)}
\otimes f^{-1}{\cal V}_{\beta }$ is unipotent since otherwise it 
is acyclic near $y$ and doesn't contribute to cohomology.

Finally we can define the combinatorial logarithmic filtered complex as:  
$$(\Psi_f ({\cal L}),{\cal W}) \colon= 
\oplus_{(\alpha ., \beta) } (\Psi^u_f 
({\cal L}^{e (\alpha .)} \otimes f^{-1}{\cal V}_{\beta }), {\cal W})) $$
 The Hodge filtration 
 $F$ extends to the logarithmic
complex and to $\Psi_f ({\cal L}_X)$.

\smallskip

\noindent { \em Theorem: Suppose ${\cal L}$ underlies a variation of polarised 
Hodge structures of weight $m$, then the complex }

{\centerline {$(\Psi_f ({\cal L}_X),{\cal W}[m], F)$}}
 
{ \em decomposes into a direct sum of intermediate extension of $VHS$; moreover we have 
$W({\cal N})= {\cal W}$ where ${\cal N}$ is the logarithm of the unipotent part of the 
monodromy $T^u$.}

\subsection {The weight filtration after M. Kashiwara and M. Saito} {\it
Local situation.} We give in this subsection  Kashiwara and Saito's constructions and 
indications on the proofs of the decomposition and the purity results for $\Psi^u_f 
({\cal L}$ [29] in order to compare the two constructions. In the reference this result 
is embedded in 
the language and theory of Hodge modules, a theory adapted for general pushforward 
results but not necessary at this stage.

\smallskip

Given $\alpha_i \in [0,1[$ for $ i \in [1,n ]$ ( or equivalently a section $\tau$ with 
value in $]0,1]$), we consider the polynomial ring 
$\C [N ]$ in one variable (resp. the field $\C [N, N^{-1} ]$) and the
module $L 
[ N ] = L {\otimes}_{\C} \C [N ]$ (resp. $L [N, N^{-1} ]$) endowed with commuting 
endomorphisms   $ (\alpha_i Id + N_i) \otimes Id$ where $N_i$ is 
nilpotent, denoted also by $\alpha_i Id + N_i$, and multiplication by $N$ denoted 
also by $N$. For each family of integers $n_i
> 0$, for $ i \in [1,n ]$, we consider the endomorphisms $A_i = \alpha_i Id + N_i- n_i N 
$ on  $L [N ]$ (resp. $L [N, N^{-1} ]$). When $\alpha_i \neq 0,\;A_i $ is invertible on 
$L [N]$ and  when $\alpha_i = 0 $, the inverse of   the endomorphisms $A_i$ are defined 
on $L [N, N^{-1} ]$  and equal to 
$$A_i^{-1} = - \Sigma_{j \geq 0} {(N_i)}^j / {(n_i N)}^{j+1}$$

\noindent where the sum is finite since $N_i$ is nilpotent for all $i$.
In 
particular $A_i$ and $A_J = {\Pi_{i \in J}} A_i \, ,\, J \subset [1, n
]$, are  
injective on $L [N ]$ so that we can deduce \par

\medskip

\noindent {\em Lemma. Given $(L, \alpha_i Id + N_i)_{i\in [1,n]} $, we
let $ 
I(\alpha.) = \{i\in [1,n]: \alpha_i = 0 \}$ 

\noindent i) The complex  $ \Omega (L[N,N^{-1} ],\; A_i = \alpha_i Id +
N_i - n_i 
N, i \in [1,n]) $ is acyclic.

\noindent ii) The following complexes are isomorphic
 $$\Omega (L [N^{-1} ], A.) \quad 
{\buildrel \sim
\over \longrightarrow}  \quad \Omega (L [N ], A.)[1 ]$$ 

\noindent iii) The complex 
$$IC (L [N ], A.) = s(Im A_{J\cap I(\alpha.)})_{J \subset [1,n ]}\simeq 0, \quad 
(Im A_{J\cap I(\alpha.)} = Im A_{J})\leqno{(33)} $$  is an acyclic sub-complex  of the  
Koszul complex $\Omega (L [N ], A.).$

\noindent iv) Let $ \Psi_J (L) :=L [N]/Im A_{J\cap I(\alpha.)}$ and ${\Psi}^{0} L = 
s({\Psi}_J (L)_{J \subset [1,n ]}, (A_i)_{i\in [1,n]})[1]$, be associated to the 
simplicial complex with differential induced by 
$A_i: {\Psi}_{J-i} (L) \to {\Psi}_J (L) $,
 then we have the following isomorphism }
$$\Omega (L[N], A.) \; {\buildrel \Pi \over {\buildrel \sim \over  
\longrightarrow}} \; \Psi^{0} L = s({\Psi}_J (L)_{J \subset [1,n ]},
A.)[1]  = s(L [N]/Im (A_{J\cap I(\alpha.)}), A.)_{J \subset [1,n ]} [1]
\leqno{(34)}$$ 

\noindent We give the statement for ${\Psi}^{\beta}(L ^{e(\alpha .)})$
in general 

\noindent {\em Proposition: Given $(\alpha .) = (\alpha_j)_{j \in 
[1,n]}$ and $ \beta $, we consider on $L ^{e(\alpha .)}$ the
endomorphisms 
 $A_i = (\tau (\alpha_i + \beta n_i)Id + N_i - n_i N) $ for $i \in [1,n]$ then 
we have the isomorphisms:}

\centerline {${\Psi}^{\beta}(L ^{e(\alpha .)}) =  \Omega (L^{e(\alpha .)}[N^{-1}], 
A.)\simeq  
\Omega (L^{e(\alpha .)}[N], A.) [1]$.}

\noindent {\em Let $ Im A_J = Im (A_{J\cap I(\tau (\alpha.+ \beta n.)}$ in $L ^{e(\alpha 
.)} $ denotes  the image of the 
composition $A_J = \Pi_{j\in J} A_j$  then this Koszul complex
is isomorphic to   }
$${\Psi}^{\beta }(L^{e(\alpha .)})\; = \;s( (L^{e(\alpha .)}[N]/Im A_{J\cap 
I(\tau (\alpha.+ \beta n.)}), A.)_{J \subset [1,n ]}[1].
\leqno{(35)}$$ 
Let $ S (\alpha .) = \{\gamma  \in {\bf C}: \forall j \in M, \alpha_j + n_j \gamma \in 
{\bf N} \}$. Only for
 $  \beta \in S(\alpha.)$ the complex is not acyclic.

\smallskip

\noindent {\em Remark i) The importance of the introduction of $ \Psi_J(L) $ is that they 
are canonically associated to the perverse sheaf $\Psi^0 (L) $, so that the 
construction of the weight filtration reduce to its construction on these vector 
spaces. It is more precise to work on these vector spaces then on the cohomology 
of the perverse sheaf, the relation being a kind of spectral sequence.}

\smallskip

\noindent {\em ii) We can give now a proof of the isomorphism (20) of the proposition in 
this paragraph. Recall
that ${\Psi}_J L = L [N ]/Im A_J$ and we have
 $Coker N/ {\Psi}_J L \simeq Coker N_J/L$  since  $A_i =
N_i - n_i N$ is equal to $N_i$ modulo $N$, so that we  
 have locally at $y \in Y^*_M$ the isomorphisms
$$ ({\bf j}_* {\cal L} / {\bf j}_{!*} {\cal L})_y\, \simeq \,{s(Coker N_J / L)}_{J 
\subset [1 , n ]} \, \simeq \, Coker N/s( {\Psi}_J L)_{J \subset [1 , n ]} \simeq Coker 
{\cal N }/ {\Psi}^u_f {\cal L}  [-1 ] $$
which establishes (20).}

\smallskip

{\em iii)  In general, the graded part of the cokernel is the primitive part
$P_k (N)$ for all  $k \geq 0 $:

 $ Gr_k ^{W(N)}(Coker N/ {\Psi}_J L ) \simeq P_k (N).$ }
 
\noindent {\em The filtration $W(N)$ on ${\Psi}_J L$  defines a filtration
by 
sub-complexes of ${\Psi}^0 L$  and  corresponds to the filtration
$W({\cal N})$ of ${\Psi}_f^u ({\cal L})$.}

Now in order to study
the weight filtration we need to consider  this complex as a perverse
sheaf in 
the corresponding abelian category. That is why we recall here basic
facts on 
this category needed to understand the construction.
 
\smallskip

 \noindent The {\it category of perverse sheaves} ${\cal L}^{\cdot}$ on
$X$ with 
respect to the  natural stratification  $Y^{\star}_M$ defined by $Y$ (
i.e such 
that for each $M \subset I$, the cohomology of  ${\cal L}^{\cdot}/Y_M^*$
is 
locally constant),  
are described locally at a point $y$ considered as the center of a
polydisc $ 
(D^*)^M$, from a topological view point, by the following combinatorial 
construction   in  [25, p 996] ( see also [16], [2]).

\smallskip

\noindent {\em The  category ${\cal P}$ of perverse sheaves  ${\cal
L}^{\cdot}$ 
on $ (D^*)^M$, with respect to its $NCD$ stratification   is  equivalent
to the  
abelian category defined as follows :\par

\noindent i) A family of vector  spaces  $L_A$ for $A \subset M$,\par 
\noindent ii) A family of morphisms \par 
\noindent $ f_{AB} \colon  L_B \to L_A \hbox { and }
h_{BA} \colon L_A \to L_B \,\hbox { for  } \,  B \subset A \subset I
\,\, \hbox 
{ such that :} $\par 
\smallskip
 \noindent $ f_{AB} \circ
f_{BC} = f_{AC}\,\,, \,\,h_{CB} \circ h_{BA} =  h_{CA}$  for  $ C 
\subset B \subset A $

\noindent $f_{AA} = h_{AA} = id\,\, , \, \,h_{A,A \cup B} \circ  f_{A
\cup B,B} 
= f_{A,A \cap B} \circ  h_{A \cap B,B} \,\,\hbox { for all} \; A,B $

\noindent and if $ A \supset B, \mid A \mid = \mid B \mid +1$, then 
 $ 1-h_{BA} f_{AB}$ is invertible. }

\smallskip
{\it Minimal extensions }

\noindent We will need the following description for $A \subset M$ of
the caegory ${\cal 
M}_A$
of the minimal extensions of a locally constant sheaf ${\cal L}$ on
$X^{\star 
}_A$: in terms of the family of vector spaces $L_B$ for $B \subset M$;
it is equivalent to $L_B = 0$ for $ A \nsubseteq B$, 
and $ f_{BA} $ is surjective and $g_{AB}$ is injective for $ A \subset B
$. We denote by ${\cal M}$ the objects isomorphic to a direct sum of
objects in $\cup_A {\cal M}_A$.

\smallskip
{\it The category ${\cal M}$ of sums of minimal extensions }

\noindent A result of Kashiwara states [25, p 997]

\noindent {\em A perverse sheaf  ${\cal L}^{\cdot}\in{\cal P}$ is a
direct 
sum of  minimal extensions (in  ${\cal M}$ ) if and only if 
$$ \forall A,  B \subset M,\;  B \subset A,\;  \quad L_A \simeq Im
f_{AB}\oplus Ker g_{BA} \leqno{(36)}$$
Moreover, it is enough to consider $ \vert A \vert = \vert B\vert + 1$. 

\noindent The above condition is equivalent to the isomorphism: $$
\oplus_{B \subset A } f_{AB} ( P_B({\cal L}^{\cdot}))\quad{\buildrel
\sim 
\over \longrightarrow}\quad L_A \leqno{(37)}$$
where  $P_B({\cal L}^{\cdot}) = \cap_{C \subsetneqq B} Ker\, g_{CB}$,
then moreover  $ 
g_{BA}:
f_{AB} ( P_B({\cal L}^{\cdot}))\rightarrow  P_B({\cal L}^{\cdot})$ is
injective 
for 
$B \subset A$.}

\medskip

{\it Description of the weight filtraion in the category of perverse
sheaves}. 

\noindent The family  $\Psi_J^{\beta} (L^{e(\alpha .)} ) = L^{e(\alpha .)} [N]/Im 
A_{J\cap I(\tau (\alpha.+ \beta n.)}$
for $J \subset M, J \neq \emptyset$ gives precisely  the description of 
$\Psi^{\beta }(L^{e(\alpha .)} )$  as a perverse sheaf, where  for $i \in J$, 
the morphisms  $f_{J(J-i)} = A_i: \Psi_{J-i}^{\beta} (L^{e(\alpha .)} 
)\rightarrow  \Psi_J^{\beta} (L^{e(\alpha .)} )$
and $g_{J(J-i)} = p_i:  \Psi_J^{\beta} (L^{e(\alpha .)} )\rightarrow 
\Psi_{J-i}^{\beta} (L^{e(\alpha .)} )$ is the canonical projection. The
product by $N$ induces on each  $ \Psi_J^{\beta} (L^{e(\alpha .)} )$ a nilpotent 
endomorphism denoted also by $N$
which  commutes with  $A_i $ and $ p_i$,  hence these morphisms are compatible 
with $W(N)$; they send $W_{r-1}(N)$ into itself (it is enough to show that  for 
$b\in \Psi_{J-i}(L^{e(\alpha .)} )$, if $ N^s (b) = 0 $ 
for $s \geq r$ , $ N^s 
(A_i(b)) = A_i (N^s (b)=0 $ (resp.for  $p_i$)). 

For each integer $r$, let $ Gr_r^{W(N)} {\Psi}_{J}^{\beta} (L^{e(\alpha .)} ), p_i^r, 
A_i^r$ denotes the
corresponding perverse graded objects, then we define 
$$ K^r_i = Ker\; p_i^r: Gr_r^{W(N)} {\Psi}_{J}^{\beta} ( L^{e(\alpha .)}) 
\rightarrow  Gr_r^{W(N)} {\Psi}_{J-i}^{\beta}  ( L^{e(\alpha .)}  ), 
\quad K^r_J:\; \cap_{i\in J} \; K^r_i \leqno{(38)}$$
in particular $\forall i \in J,  K^r_J \subset Ker N_i \subset
Gr_r^{W(N)} {\Psi}_{J}^{\beta} ( L^{e(\alpha .)})  $.
The aim of the next part is to deduce the decomposition property (37) via the 
proof of (36) in presence of a polarised Hodge filtration.

First we give a global setting of the problem.  

\be{The global weighted complex of nearby cycles  $({\Psi}_f
(\cal L), W({\cal N})$}{3.II.1} \
In this subsection, we define the weight filtration abstractly without going back to an 
explicit formula as used on the combinatorial logarithmic complex.
 The filtration $W(N)$ on each 
${\Psi}^{\beta } (L^{e(\alpha .)})$  defines a filtration by
sub-complexes on 
$\oplus_{\alpha.} {\Psi}^{\beta }  L^{e(\alpha.)}$ and corresponds via
(32) to 
the filtration $W({\cal N})$ of ${\Psi}_f ({\cal L})_y$ in the abelian
category 
of  perverse sheaves where ${\cal N} = -1/2i\pi Log T^u $. 

\noindent Consider  Deligne's extension  $({\cal L}^{e (\alpha
.)}\otimes 
f^{-1}{\cal V}_{\beta })_X$ and the associated logarithmic complex 

\centerline {$ i_Y^{\s} s({\Omega}^*_X (Log Y) \otimes ({\cal L}^{e (\alpha.)}\otimes 
f^{-1}{\cal V}_{\beta })_{X}[p], \eta )_{p \geq 0}$.}

There exists a global acyclic sub-complex $IC(({\cal L}^{e (\alpha.)}\otimes f^{-1}{\cal 
V}_{\beta })_X[N])$
inducing at each fiber at the point $y$ the complex 
$IC(L^{e (\alpha .)}[N], A_i = (\tau (\beta n_i +\alpha_i) + N_i - n_i N))$. We define 
the weight filtration $ W({\cal N})$ on the quotient complex
$$\Psi^u_f (({\cal L}^{e (\alpha .)} \otimes f^{-1}{\cal V}_{\beta})_X)\colon= 
(s({\Omega}^*_X (Log Y) \otimes ({\cal L}^{e (\alpha .)}\otimes f^{-1}{\cal 
V}_{\beta })_{X}[p], \eta )_{p \geq 0} {\hbox {\large /}}
IC({\cal L}^{e (\alpha .)}\otimes f^{-1}{\cal V}_{\beta })_X[N]))[1]$$
as the filtration inducing $W(N)$ at each fiber 
${\Psi}^{\beta } (L^{e(\alpha .)})$   at points of $Y$.
Finally we can define the logarithmic filtered complex as:
$$(\Psi_f ({\cal L}_X),W({\cal N}))\colon=\oplus_{(\alpha ., \beta) }(\Psi^u_f (({\cal 
L}^{e (\alpha .)} \otimes f^{-1}{\cal V}_{\beta }))_X,W({\cal N})) 
\leqno{(39)}$$
The filtration $F$ extends to $(\Psi_f ({\cal L}_X)$ via its extension to
${\Omega}^*_X (Log Y) \otimes {\cal L}_X [N][1]$ by the formula 
 $$ F^p(i_Y^{\s} s({\Omega}^*_X (Log Y) \otimes {\cal L}_X[i])_{i\geq 0})[1] = i_Y^{\s} 
s(F^{p+i+1}(\Omega^*_X (Log Y) \otimes {\cal L}_X ) [i], \eta)_{i\geq 0}[1]$$
$$ = i_Y^{\s} s[(s(\Omega^q_X (Log Y) \otimes F^{p+i+1-q}{\cal L}_X )_{q \geq0})[i], 
\eta)]_{i\geq 0}[1].$$

\noindent { \em Theorem: Suppose $\cal L$ underlies a variation of polarised 
Hodge structures of weight $m$, then graded part of the weight on the complex (39)
with the filtration $F$ defined above
$$(\Psi_f ({\cal L}_X),W({\cal N})[m], F)\leqno{(40))}$$ 
decomposes into intermediate extensions of $VHS$.}

\smallskip

The proof of this theorem reduces by definition to show that
$(Gr^{W({\cal N})[m]}_r,F) $ decomposes.
This can be checked locally via (36) and (37). That is we need to use the following 
decomposition theorem based  on results due to Kashiwara [25] in  characteristic zero and 
 proved in [2] in the language of purity in positive characteristic. 

\smallskip

{\em Theorem (decomposition)(Kashiwara - Saito): For each onteger $a, \; Gr^{W({\cal 
N})}_a (\Psi_f ({\cal L}^{e(\alpha.)}))_y \simeq Gr^{W(N)}_a 
\Psi^{\beta }(L^{e(\alpha .)})$ is isomorphic to a direct sum  of fibers at $y$ 
of various  intermediate extension of  variations of polarised Hodge structures. 
Precisely 
$$  Gr^{W(N)}_a \Psi^{\beta }(L^{e(\alpha .)})\simeq  \oplus_{J\subset M} IC 
((K^a_J(L^{e(\alpha .)}), N_i, i\in M-J)\leqno{(41))}$$
where $K^a_J$, defined by $(38)$, is a pure Hodge structure of weight
$a + m$ with the induced Hogde filtration $F$.}

\smallskip

i) Elements of Kashiwara's proof [30, prop. 3.19, and Appendix]. We will
write $L$ for $L^{e(\alpha .)}$ 
and associate to $ (L, F, P, N_{i}, i\in M = [1,n] )$ the module $L[N]$ where $N$ is a 
polynomial variable, endowed with two filtrations as follows. Consider $W(L) = 
W(\Sigma_{i \in M} N_i) [m ]$ and $F$ on $L$, then define 
$$W_k(L[N ]) = {\Sigma}_j W_{k+2j} L \otimes N^j \, , \, F^p (L[N]) = {\Sigma}_j F^{p+j} 
L \otimes N^j \leqno{(42)}$$
Since the endomorphisms $A_i = N_i - n_i N$ shift $W$ by $-2$ and $F$ by $-1$, the two  
filtrations induce a $MHS$ 
 on the cokernel $\Psi_J = L[N]/Im A_J$ for $J \subset M $.
 We have  an isomorphism compatible with  the filtrations
$$ ({\oplus}_{j \leq l-1} L \otimes N^j, W, 
F) \quad \simeq \quad ({\Psi}_J L, W, F) \leqno{(43)}$$
 obtained via the composition of the natural embedding in $ L[N]$ with the projection on 
${\Psi}_J L$, where $W$ anf $F$ are defined on the left term as in the formula
(42) above.
 In fact the relation : $N^l = \Sigma_{j \in [1,l]}(-1)^{j+1} \sigma_j (
\frac{N_i}{n_i}, i \in J)\otimes N^{l-j}$ where $l = \mid J \mid$
and ${\sigma}_j$ is the $j^{th} $ elementary symmetric function of 
$\frac{N_i}{n_i}, i \in J$, on the quotient of the right term
leads to the definition of the action of
$N$  on the left term by the formula :

\centerline { $N (a \otimes N^{l-1}) = 
\Sigma_{1 \leq j \leq l}{(-1)}^{j+1} {\sigma}_j ((N_i/n_i),i \in J) (a) \otimes 
N^{l-j}$.}

 In order to define a polarisation we introduce 
a product $P_J$ on ${\Psi}_J (L)$ as follows

\centerline { $P_J (a N^i, b N^j) = P(a, {(-1)}^i \quad Res\; (A_J^{-1}(b \otimes 
N^{i+j}))) $}

\noindent where $A^{-1}_J$ is defined on $L [N,N^{-1}]$, $N$ is
considered as a variable $x$ and the residue $\, Res\,$ is equal to the
coefficient of $1/N$ in the  fraction in $N$. This formula shows
directly that the product is well 
defined on $Coker A_J$ ; in fact, $P_J (a N^i, A_J (c) = P(a, {(-1)}^i \quad \hbox{Res} 
(c \otimes N^i) = 0 $ since the residue is zero.  Using an explicit 
expression of $A^{-1}_i$, we find $P_J(a N^p, b N^q)= (-1)^p P(a,\Sigma_{a_i} 
(\Pi_i (N_i^{a_i}(b)/ n_i^{a_i+1}))$ where $a_i \geq 0$ and
$\Sigma_i a_i = i+j-l+1$. In particular
 
\centerline {$ P_J (a , b N^r) =
(1/\Pi_i n_i) P(a,b) \, {\hbox { if}}\, r = l-1,$ and zero
otherwise}

\centerline {$ P_J (a N^i, b N^j) ={(-1)}^i P_J(a, b N^{i+j})$. }

 In $[30]$, the following result is attributed to Kashiwara

\smallskip

\noindent {\em Theorem : With the previous notations, namely $W$ and $F$ 
$$\Psi_J(L) = (L[N]/Im A_J, N_1,\ldots, N_n, N; W,F, P_J)
\leqno{(44)}$$
 underlies  a polarised nilpotent orbit of weight $m+1 - \vert J \vert,$ that is:  the 
weight filtration $ W(N + \Sigma_{i\in M} N_i)[m+1 - \vert J \vert] = W $
underlies the weight of a $MHS$ on $\Psi_J(L)$ 
with the Hodge filtration $F$. }

\smallskip

ii)(see M. Saito [29,5.2.15, 5.2.14], [30,3.20.4]). The induced morphisms $N, N_i$ and 
$A_i $ shift $W$ by $-2$ and $F$ by $-1$.
 Since  $W(N + \Sigma_{i\in J} N_i)$ is  the  weight filtration of the endomorphism 
$\Sigma_{i\in J} N_i$ relative to $W(N)$  that is for all $ \Psi_J(L)$ :

\smallskip

\centerline {$Gr^{W(N + \Sigma_{i\in J} N_i)}_{j+r} Gr^{W(N)}_r {\buildrel (\Sigma_{i\in 
J} N_i)^j \over  \longrightarrow } Gr^{W(N + \Sigma_{i\in J} N_i)}_{r-j} Gr^{W(N)}_r $}

we have: 

\centerline {$Gr^{W(N + \Sigma_{i\in J} N_i)}_{j+r} Gr^{W(N)}_r \simeq Gr^{W(\Sigma_{i 
\in 
J} N_i)}_j Gr^{W(N)}_r \simeq Gr^{W(\Sigma_{i \in 
J} A_i)}_j Gr^{W(N)}_r$}

\smallskip

Now we may consider the orbit with only two endomorphisms $ (\Psi_J(L), N_i, N, F = 
F(N_i))$ ($F$ is the limit along the axis $Y_i$), then we deduce  
commutative diagrams for $j$ varying in an interval of {\bf Z} symmetric with 
center $0$ with at left $HS$
of weight $n+j$ where $ n = m + a - \vert J \vert $ and to the right $n + j - 1$

$$\matrix{Gr_{j}^{W(N_i)} Gr_a^{W(N)} {\Psi}_{J-i} L&{\buildrel \rm A_i
\over \longrightarrow}&Gr_{j-1}^{W(N_i)} Gr_a^{W(N)} {\Psi}_{J} L (-1)\cr \qquad 
\downarrow{N_i}&\swarrow{p_i}&\downarrow{N_i}\cr
Gr_{j-2}^{W(N_i)} Gr_a^{W(N)} {\Psi}_{J-i} L(-1) &{\buildrel \rm A_i \over 
\longrightarrow}&Gr_{j-3}^{W(N_i)} Gr_a^{W(N)} {\Psi}_{J}  L (-2 )\cr} \leqno{(45)}$$
\noindent  moreover we have:
$\, P_J(A_i u,v) \,=\,P_{J-i}(u,p_i v) \,$ for all $\,u \in {\Psi}_{J-i}
L \,$ and $\,v \in {\Psi}_{J} L \,$ \par  
\noindent  In this situation, a result of M. Saito [MI,5.2.15]  applies and shows

\smallskip

\noindent {\em Proposition :
 For all $J \subset I$ and $i \in J$,
consider the  morphisms

\smallskip

\centerline {$Gr_a^{W(N)} {\Psi}_{J-i} L {\buildrel A_i \over \longrightarrow} \, 
Gr_a^{W(N)} {\Psi}_{J} L (-1)\,{\buildrel p_i \over \longrightarrow}\, Gr_a^{W(N)} 
{\Psi}_{J-i} L (-1)$}

\noindent then we have  a  decomposition  

$$Gr_a^{W(N)} {\Psi}_{J} L \quad \simeq \quad Im A_i  \oplus ker \, p_i 
\leqno{(46)}$$

\noindent compatible with the primitive  decomposition . In particular, $p_i$ 
induces an isomorphism of $Im A_i$ in $Gr_a^{W(N)} {\Psi}_{J} L$ onto $Im A_i$ 
in $Gr_a^{W(N)} {\Psi}_{J-i} L$.} \par

\smallskip

\noindent The result is deduced from the sequence in the proposition by taking its graded 
version  $Gr_{j-1}^{W(N_i)}
$ for various $j$ as in $(45)$ and using the polarisation of $HS$ induced on, to prove 
for each $j$,

\centerline {$ Im A_i \oplus Ker p_i  \simeq Gr_{j-1}^{W(N_i)} Gr_a^{W(N)}{\Psi}_{J} L 
(-1)$}

 hence, since $A_i$ and $p_i$ are compatible with the $MHS$
of weight $W(N_i)$, we get:

\noindent $ Im A_i \oplus Ker p_i  \simeq  Gr_a^{W(N)} {\Psi}_{J} L (-1)$. 

\smallskip

Now, to finish the proof of the decomposition theorem, it remains  to show that $ K^a_J$ 
is pure and polarised in two steps:

\smallskip

 \noindent {\em Lemma: i)$ K^a_J \subset W_0(  \Sigma _{i\in J} N_i) Gr_a^{W(N)} 
{\Psi}_{J} L $.

\noindent ii) $ K^a_J \cap (W_{ - 1}( \Sigma_{i\in J} N_i)
Gr_a^{W(N)} {\Psi}_{J} L ) = 0$.}

\smallskip

\noindent Proof. i) the assertion (i) follows from the relation: $ Ker
(\Sigma_{i\in J} N_i) \subset W_0(\Sigma_{i\in J} N_i) 
Gr^{W(N)}_a {\Psi}_{J} L$.

\noindent ii) Suppose $x \in W_{-s}(\Sigma_{i \in J} N_i) Gr^{W(N)}_a {\Psi}_{J} L \cap 
K^a_J $ where $ - s \leq -1 $, then there exists \\
 $y \in W_{s}(\Sigma_{i \in J} N_i)
 Gr^{W(N)}_a {\Psi}_{J} L \; = W_{s}(\Sigma_{i \in J} A_i)
 Gr^{W(N)}_a {\Psi}_{J} L $ such that  $ x = (\Sigma_{i \in J} A_i)^s (y) $ (by 
surjectivity of $\Sigma_{i \in J} A_i$ on negative weights ) then  for each $i$, we have  
 $(N_i)^s 
(y)\in  Im A_i \; mod\; W_{- s-1}$, hence   $ x = \sum_{i} N_i^s (y)$ is in 
$ (\cap_i ker p_i) \cap \sum_{i} Im A_i  = 0 \; mod\;  W_{s-1}$, that is $x
\in W_{-s-1}(\Sigma_{i \in J} N_i) Gr^{W(N)}_a {\Psi}_{J} L $. We deduce (ii) by a 
descending  inductive argument on $- s$.

\smallskip

We deduce from the lemma that $K^a_J $ is pure of weight $a$
which ends the proof of the theorem.

\smallskip

\subsection { Example : Rank one local system $ {\cal L} $ on $X-Y$.}
We apply the above  theory to remove the base change in 
Steenbrink's work.
 In this case the monodromy of $ {\cal L} $ around components of $Y$ is of the following 
form: $\forall i \in I, T_i = \alpha_i Id, N_i = 0 $. For $\beta \in \Q\cap [0,1[$, let 
${\cal V}_{\beta }$ denotes the rank one local 
system on the punctured disc with monodromy $ e^{-2i\pi \beta }$, $ {\cal S}\;\colon= 
\;{\cal L} \otimes f^{-1}{\cal V}_{\beta } $ and ${\cal S}_X$ its Deligne's extension;
then we can define the weight filtration $ W({\cal N})$ explicitly on the complex
$$\Psi^u_f (({\cal S}_X)\colon= i_Y^{\s} s({\Omega}^*_X (Log Y) \otimes {\cal S}_X ) [p], 
\eta )_{p \geq 0}[1] $$
 First we define $W$ on $\Omega^*_X (Log Y) \otimes {\cal S}_X $. Let $ I(\beta):\{i\in 
I: \alpha_i + \beta n_i \in \Z\}$, $Y(\beta) = \cup_{i\in I(\beta)}Y_i$, $C(\beta)= 
\cup_{i\in I-I(\beta)}Y_i$, notice that ${\cal S}_X$ is locally trivial along $Y(\beta) - 
C(\beta) $ and its logarithmic complex is acyclic along $C(\beta) $,
that is for $j: X-Y \rightarrow X$, we have:  $({\bf j}_* {\cal S})_{/X-Y(\beta)})= ({\bf 
j}_! {\cal S})_{/X-Y(\beta)}$. 

\noindent We write  $\Omega^*_X (Log Y)$ as $\Omega^*_X (Log Y(\beta) \otimes \Omega^*_X 
( Log \,C(\beta))$ and extend the logarithmic weight filtration $W^{Y(\beta)}$ along 
$Y(\beta)$ to the whole complex by 

\noindent $ W \; \colon= \;[W^{Y(\beta)}(\Omega^*_X (Log Y(\beta))]\otimes \Omega^*_X 
(Log C(\beta))\otimes {\cal S}_X $. Then we have 

\noindent $Gr^W_i (\Omega^*_X (Log Y)\otimes {\cal S}_X ) \simeq 
(Gr_i^{W^{Y(\beta)}}(\Omega^*_X (Log Y(\beta))) \otimes \Omega^*_X (Log C(\beta))\otimes 
{\cal S}_X \simeq$ 

\noindent $\oplus_{J \subset I(\beta ), \mid J\mid = i} \Omega^*_{Y_J}
(Log(C(\beta)\cap Y_J) )[-i] \otimes {\cal S}_{Y_J} $

\noindent with the differential of the induced connection on ${\cal S}_{Y_J} $.

\noindent Locally, let $L$ denotes the general fiber of ${\cal S}$, then the fiber 
of the logarithmic complex at a point $y \in Y^*_M$ is isomorphic to the Koszul complex
$(\Omega(L, \tau (\alpha_i + \beta n_i)Id, i \in M )$ where $L_J$ corresponds to
$L \otimes \wedge_{i \in J} \frac {dz_i}{z_i}$.
 For each $i\in M$, let $(i,\beta ) = \{ j\in I-I(\beta ):
Y_j\cap Y_i\neq \emptyset\}$,  then the  fiber of $Gr^W_i$ is:

\centerline {$Gr^W_i \simeq  \oplus_{i \in I(\beta )} (\Omega(L, \tau (\alpha_j + \beta 
n_j)Id, j \in (i,\beta ))$ }

\noindent This weight filtration and  the Hodge filtration extend to $s({\Omega}^*_X (Log 
Y) \otimes {\cal S}_X [p], \eta )_{p \geq }[1] $ by the formula 

\centerline {$W_i = \oplus_{p \geq 0} W_{i + 2p + 1}, \;\; F^i = \oplus_{p \geq 0} F^{i + 
p +1}$. }

\noindent Here  
 the fiber of the double complex at $y$ is isomorphic to the Koszul complex
$(\Omega(L[N], A_i = \tau (\alpha_i + \beta n_i)Id - n_i N)$ where $A_i$ is an 
isomorphism whenever $\tau (\alpha_i + \beta n_i)\neq 0$, that is $i \in I - I(\beta )$.
                 
\noindent Now we introduce the acyclic subcomplex ${\cal K} \simeq \oplus_{p \geq 0} 
W_p[p+1] $, whose fiber at $y$ in
$(\Omega(L[N], A_i = \tau (\alpha_i + \beta n_i)Id - n_i N, i\in M)$ 
is given as $s(K_J)_{J\subset M}$ where $K_J = Im A_J = N^r L[N]$
whenever $\mid J \cap I(\beta ) = J(\beta )\mid = r \; (A_J = \Pi_{i \in J}A_i)$. 

\noindent The quotient complex with induced filtration:
$$\Psi^{\beta }_f ({\cal L}) \cong ( i_Y^{\s}[s({\Omega}^*_X (Log Y)\otimes {\cal 
S}_X [p], \eta )_{p \geq 0}[1]]/{\cal K}, W , F) $$

\noindent is the bifiltered  complex computing $\Psi^{\beta }_f ({\cal L})$.
The quotient complex has fiber at $y$: 
$$s(\Psi _J(L), A_i, i\in M)_{J\subset M} , 
\Psi_J(L) = (L[N] /ImA_J ) \simeq {\oplus_{p \in [0, r-1]}L N^p}.$$
\noindent Let ${\cal K}_p$ denotes the $pth $ column of ${\cal K}$, then the graded 
object with respect to $W$ is:

\smallskip
\centerline {$ Gr^W_i (\Psi^{\beta }_f ({\cal L}) \simeq
i_Y^{\s} s(Gr^{W}_{i+2p+1}[(\Omega^*_X 
(Log Y)\otimes {\cal S}_X )[p]]/{\cal K}_p)_{p \geq 0}[1]   \simeq $}

\smallskip

\centerline {$ s(\oplus_{J \subset I(\beta ), \mid J\mid = i+ 2p +1} \Omega^*_{Y^J}(Log 
(C(\beta)\cap Y^J) )\otimes {\cal S}_{Y^J}[-i - p ])_{p\geq 0,p\geq -i, p \leq n= dim \; 
X}   $}

\smallskip

\noindent {\em Proposition: $_WE_1^{p,q}((\Psi^{\beta }_f ({\cal L}) ) = \oplus_{i\geq
0,i\geq p} \oplus_{J\subset \beta } H^{2p+q-2i}_c ( Y_J - C(\beta ))^{-p+2i+1}, {\cal 
S}_/)$

\noindent where ${\cal S}_/$ is a local system on $Y_J$ deduced from ${\cal S}$ 
as restriction of Deligne's extension ${\cal S}_X$ and has zero restriction to $C(\beta 
)$. }

\noindent The weight spectral sequence is: $_WE_1^{p,q}((\Psi^{\beta }_f ({\cal L}) ):= $

\noindent $ H^{p+q}(X, Gr^W_{-p}((\Psi^{\beta }_f ({\cal L}) ) = \oplus_{i\geq 0,i\geq p} 
H^{p+q-i} (X, Gr^{W}_{-p+2i+1}
(i_Y^{\s}[(\Omega^*_X (Log Y)\otimes {\cal S}_X )[i]]/{\cal K}_i) [1]) = $

\noindent $ \oplus_{i\geq 0,i\geq p} H^{p+q+1} (\tilde Y,
\Omega^*_{\tilde Y(\beta ) ^{-p+2i+1}} 
(Log (C(\beta)\cap \tilde Y(\beta )^{-p+2i+1})
 \otimes {\cal S}_{\tilde Y(\beta ) ^{-p+2i+1}}  [p-2i-1]) = $

\noindent  $ \oplus_{i\geq 0,i\geq p} H^{2p+q-2i}_c (\tilde Y(\beta )-C(\beta )) 
^{-p+2i+1}, {\cal S}_/)$
where ${\cal S}_/ \simeq {\cal S}_{\tilde Y(\beta ) ^{-p+2i+1}}  $ is induced by ${\cal 
S}_X$ and has zero restriction to $C(\beta )$. 

\section{Variation of Mixed Hodge structures}
Let $(L, W^0)$ be a filtered  object in an abelian category
and $N$ a nilpotent endomorphism of $(L,W^0)$. Deligne [10, (6.1.13)]
introduced the notion of relative weight filtration
$W$ of $N$  with respect to $W^0$ on $L$ and showed that if it exists,
it is the unique filtration satisfying for all $a\in \Z, b\in {\bf N}$

\smallskip

\centerline {$ N W_a \subset W_{a-2} \; \hbox{and}\; N^b: Gr ^{W}_{a+b}Gr^{W^0}_a L 
\simeq Gr ^{W}_{a-b}Gr^{W^0}_a L $}

\smallskip

\noindent A variation of mixed  Hodge structures $ (VMHS ):\; (L,W^0, F)$ is called 
  good  [10,(1.8.15)] if there exists a relative weight filtration
$W$ for the action of  the logarithm of the monodromy $N$.
 We showed in three notes developed in [13], the existence of a limit
relative weight filtration $W$ for geometric $VMHS$
inducing on $Gr^{W^0} L$ the limit of the $VHS$ on 
$(Gr^{W^0} L, F)$. Steenbrink and Zucker called an axiomatic  $VMHS$
admissible if it is good and satisfy a set of  properties
all satisfied by the geometric case [25], [34].
In this section we show that the definition of the weight 
filtration extends to this case without major difficulties
using only the ingredients of proofs  already introduced
in the previous cases.

\be { Good $VMHS$.}{4.1} \

Let $ V_X = ({\cal L}, ({\cal L}^{\Q}, W^0 ),
 ({\cal L}^{\C}, W^0\otimes\C, F ))$ be a 
unipotent $VMHS$ 
on $X-Y$ and ${\cal L}_X$ its canonical extension, then $W^0$ (finite) 
extends to a filtration by subbundles. We  say that $ V_X$ is
good if the following properties are satisfied :

\noindent i) the filtration $F$ extends to ${\cal L}_X$  as a filtration by sub-bundles

\noindent ii) for all $J \subset I$, the relative filtration $ M_J \colon = W(N_J, 
W^0_{Y_J^*})$ exists  ( it is a filtration by sub-local systems).

\noindent iii) The limit filtrations $(M_J, F)$ define a $VMHS :\; V_{Y_J^*}$ on $Y_J^*$; 
moreover 
$W^0_{Y_J^*}$ is a filtration by sub-$VMHS$ such that the induced $VMHS$  
on $Gr^{W^0} {\cal L}^{\Q}_{Y_J^*}$ coincides with the limit of the $VHS$
on $Gr^{W^0} {\cal L}^{\Q}_{X^*}$ ($V_{Y_J^*}$ is called the limit of 
$V_{X^*}$ on $ Y_J^*$).

\noindent iv) Compatibility: let $K,J \subset I $, then the $VMHS \;\; V_{Y_J^*}$
(iii) satisfy (ii) on $Y^*_J$ and its limit $VMHS$ on $ Y_{J \cup K}^*$
 coincides with the $VMHS:\; V_{Y_{J \cup K}^*}$ limit of the $VMHS: \;
 V_{X^*}$, that is to say
 $$ W(N_{J \cup K}, W^0_{Y_{J \cup K}^*}) =  W(N_K, W(N_J, W^0_{Y_J^*})).$$
The last property is to be understood at each point $y \in Y_{J \cup K}^* $ where $ 
W(N_J, W^0_{Y_J^*})$ extends, moreover the above properties are not independant. In a  
study 
[25], Kashiwara deduce the properties (ii) to (iv) from 
the existence of $M_i$ for $i \in I$.
 
\noindent The $VMHS$ is said to be  graded polarised if for all $r, \;Gr^{W^0}_r {\cal 
L}$ is polarised.

\be { Local definition of the weight filtration.}{4.2} \

For each
subset $s_{\lambda} \subset I$ such that $Y_{s_{\lambda} } \neq 
\emptyset$ we write $ {\cal W}^{s_{\lambda}} = W (\Sigma_{i\in 
{s_{\lambda}}} {\cal N}_i, W^0)$ for the filtration by subbundles defined by the 
nilpotent endomorphisms of the restriction 
${\cal L}_{Y_{s_{\lambda}}}$ of $ {\cal L}_X$ to $Y_{s_{\lambda}}^*$ (which exists by 
hypothesis), inducing 
also by restriction and for each subset
 $K \supset  s_{\lambda}$ a
filtration on $ {\cal L}_{Y_K}$. We define the weight filtration 
${\cal W}' (s.)$ for $s. \in S(I)$, locally near each point $y \in Y^*_M $, as in the 
case of $VHS$, in terms of a set of $n$ coordinates $y_i, i\in [1,n]$ where we
identify $M$ with $[1,p]$ on an open set $U_y\simeq D^{\mid M\mid}\times D^{n-p}$ 
containing $y$ and  a section $f^{s.} = \Sigma_{J\subset M,J'\cap  M = \emptyset} 
f^{s.}_{J,J'} 
\frac {dy_J}{y_J}\wedge dy_{J'}$ of 
${\Omega}_{X_{s.}}^* \otimes {\cal L}_{X_{s.}}$

$$ f^{s.}\in {\cal W}'_r ({\Omega}_{X_{s.}}^* \otimes {\cal L}_{X_{s.}})_{/U_y}
\Leftrightarrow \forall J, N \subset M,\, f^{s.}_{J,J'}/Y_N \cap U_y \in {\cap \atop 
{s_{\lambda}\in s.,s_{\lambda} \subset N}} {\cal W}_{a_{\lambda}(J,r)}^{s_{\lambda}} 
{\cal L}_{(Y_N\cap U_y)} $$
\noindent It is a filtration by subcomplexes of analytic subsheaves globally defined on 
$X$. 
We deduce two filtrations  ${\cal W}'$ and ${\cal W}$  of 
${\Omega}^* {\cal L}$ as follows:
 $${\cal W}' = s({\cal W}' (s.))_{s. \in S(I)}, \;\; {\cal W}_r (s.) = {\cal W}'_r (s.) 
\cap W^0_r, \;\; {\cal W} = s({\cal W}(s.))_{s. \in S(I)}$$

moreover the filtration $W^0$ extends as a constant 
filtration for all $s.\in S(I)$ and the previous definition of the Hodge filtration $F$ 
remains unchanged. The filtrations
 ${\cal W}'$ and ${\cal W}$ will have different applications,
the first leads to the filtration ${ W}({\cal N}) $
on the nearby-cocycles and the second defines a $MHS$ on 
$X-Y$ when $X$ is proper. 

Let $b \in \Z $ and $K \subset M $; when the point  $y $ is in $ Y^*_K \subset Y^*_M$,
the previous study apply to the nilpotent orbit $Gr^{W^0}_{b}L$ defined at
the point $y$ so that we can conclude that for $a > b$ (resp.  $ a < b $) the complex 
$C^K_{a} (Gr^{W^0}_{b}L)$ is concentrated in degree $\mid K \mid$ (resp. $\mid K\mid - 1$ 
 and acyclic for $ a = b$. When the point 
$ y $ vary in $ Y^*_K$, the cohomology       $ H^{\mid K \mid}(C^K_{a} Gr^{W^0}_{b} L)$ 
is of weight $ a - \mid K \mid $ induced by $ W^K_{a-\mid K \mid} (N, W^0 L)$ and $F$  ( 
resp.  $H^{\mid K\mid - 1}(C^K_{a} (Gr^{W^0}_{b}L)))$
of weight $ a + \mid K \mid $),
defines a local system ${\cal L}^K_{a, b} $ on $ Y^*_K$ which underlies  a $VHS$ induced 
by $F$.
We define as well the local system ${\cal L}^K_a$ of 
general fiber the cohomology of  $C^K_{a} (W^0_{a-1}L)$ equal to   $ H^{\mid K 
\mid}(C^K_{a} W^0_{a-1} L)$ is of weight $ a - \mid K \mid $ induced by $ W^K_{a-\mid K 
\mid} (N, W^0 L)$ and $F$, underlies a   polarised $VHS$  induced by $F$ graded by $W^0$.

\smallskip
\noindent { \em Theorem. Let  $(\cal L, W^0, F)$ be
a unipotent graded polarised good variation of mixed
Hodge structures  on X-Y; then the complex 

\centerline {$ ({\Omega}^* {\cal L}, {\cal W}', {\cal W}, W^0, F)$}
 
\noindent  with the filtrations ${\cal W}', \; {\cal W}, \; W^0$ and $F$ defined above 
satisfy 
the following decomposition and purity properties

\smallskip

i) Purity:  For all subset $K \subset I$ and all integers $a > b $ (resp. $a < b $), the 
$VHS: \;({\cal L}^K_{a, b}, W, F) $ 
 of general fiber  $ H^{\mid K\mid} (C^K_{a} (Gr^{W^0}_{b}L))$ ( resp.  $H^{\mid K\mid - 
1}(C^K_{a} (Gr^{W^0}_{b}L)))$ on $Y^*_K $,
 as well the $VHS: \; ({\cal L}^K_a, W, F)$ of 
general fiber  $C^K_{a} (W^0_{a-1}L)$, graded  polarised by $W^0$. 

\smallskip
 
ii) The complex $ Gr^{{\cal W}'}_{a}Gr^{W^0}_{a}{\Omega}^*{\cal L}$ is acyclic, hence

\centerline { $ Gr^{{\cal W}}_{a}{\Omega}^*{\cal L}\simeq Gr^{{\cal W}'}_{a}W^0_{a - 
1}{\Omega}^*{\cal L} \oplus Gr^{W^0}_a {\cal W}'_{a-1}{\Omega}^*{\cal L}$ }

\smallskip

iii)Decomposition : We have the following decomposition 
into intermediate extensions on $Y_K$ of $VHS$ ${\cal L}^K_{a, b} $ ( resp. ${\cal 
L}^K_a$) 

\centerline {$( Gr^{{\cal W}'}_{a}Gr^{W^0}_{b}{\Omega}^* {\cal L}, F) \simeq  \oplus_{K 
\neq \emptyset, K\subset I} j_{!*}^K ({\cal L}^K_{a,b}
[-\mid K \mid], W[ 2 \mid K \mid], F[- \mid K \mid])\;\;{\hbox {if }} a > b, \;\; 
(j^K:Y^{*}_K \rightarrow Y_K), $}

\centerline  {$( Gr^{{\cal W}'}_{a}Gr^{W^0}_{b}{\Omega}^* {\cal L}, F)  \simeq  \oplus_{K 
\neq \emptyset, K\subset I} j_{!*}^K ({\cal L}^K_{a,b} 
[1 - \mid K\mid], W[-1] , F)\;\; {\hbox {if}} a < b $,}
 
\centerline { $ Gr^{{\cal W}}_{a}{\Omega}^* {\cal L} \simeq  \oplus_{K \subset I} 
j_{!*}^K ({\cal L}^K_{a} [-\mid K\mid ], W[2 \mid K \mid], F [- \mid K \mid])$.}

\smallskip

iv) When $X$ is proper, the filtrations ${\cal W}$ anf $ F$  define a $MHS$ on the 
cohomology of $X-Y$ with value in ${\cal L}$  and the filtration $W^0$
induces a filtration by $sub-MHS$.}

\smallskip

Proof. i) Locally we reduce the problem to the study for 
$K \neq \emptyset , K \subset M \subset I $ of
$C^{K}_a (Gr^{W^0}_b L)$ (resp. $C^{KM}_a (Gr^{W^0}_b L)$) previously seen for the 
polarised nilpotent orbit $Gr^{W^0}_b L$, while the fiber of 
$j_{!*}^K {\cal L}^K_{a}$ at $y$ is the complex $C^{KM}_a (W^0_{a-1} L)$. The complex 
$C^{K}_a (W^0_{a-1} L)$ for a nilpotent orbit $L$ defined at $y$ in $Y^*_K$ is a $VHS$ by 
successive extensions of $Gr^{W^0}_i L$ for $i < a$.

\smallskip

ii) The acyclicity reduces to the case of the $VHS: \;  Gr^{{\cal 
W}'}_{a}{\Omega}^*Gr^{W^0}_{a}{\cal L} \simeq 0$,
while the direct sum is similar to the case studied in ($\S 2. II$). Notice that: 
$Gr^{{\cal W}'}_{a}W^0_{a-1}{\Omega}^*{\cal L} = Gr^{{\cal W}'}_{a}{\Omega}^* 
W^0_{a-1}{\cal L} $.

\smallskip
iii) The assertions for $ Gr^{W^0}_b L$ follow from the case of $VHS$, while the 
assertion for ${\cal W}$ follows from (ii)
whose right term appears in the case of  $K = \emptyset $ and  $Gr^{{\cal W}'}_{a} 
{W^0}_{a-1} L$ appears for  $K  \neq \emptyset $ and can be checked  as successive  
extensions of $Gr^{W^0}_i L$ for $i < a $ using the count of weight as in the case of 
$VHS$ and following remark 
  
\smallskip 
{\em Remark: We use in the proof the fact that a bifiltered complex  
$(K, W^0, F)$ with $(K,W^0)$ defined over $\R$ and a finite increasing filtration $W^0$ 
such that  $(Gr^{W^0}_r K, F)$ is a Hodge 
complex of weight $a$ for all $r$, then $(K,F)$ is a Hodge 
complex of weight $a$ and $W^0$ induces on cohomolgy 
a filtration by $sub-HS$.}

\smallskip
iv) When $X$ is proper, the decomposition and purity results prove that 
 the complex $ ({\Omega}^* {\cal L}, {\cal W}', W^0, F )$ is a filtered mixed Hodge 
complex  according to the terminology of [13].

\be { Nearby-cocycles $\Psi^u_f {\cal L}$}{4.3} \

The constructions for  $\Psi^u_f {\cal L}$ are similar to the case of $VHS$. 
We define the filtrations

\centerline { $ {\cal W}_r' ({\Psi}^{u}_f {\cal L})_{ S(M)} = i_Y^{\s}s ({\cal 
W}_{r+2p-1}' {\Omega}^* {\cal L} [p], \eta)_{p \leq 0},$}

\centerline { $ W_r^0 ({\Psi}^{u}_f {\cal L})_{ S(M)} = i_Y^{\s}s ( {\Omega}^* 
(W^0_r{\cal L}) [p], \eta)_{p \leq 0},$}

\centerline { $F^r ({\Psi}^{u}_f {\cal L})_{ S(M)} = i_Y^{\s} s ( F^{r+p} {\Omega}^* 
{\cal L}[p], \eta)_{p \leq 0} $}

the action of the logarithm of the monodromy $\nu$ is defined similarly,
 then the relation for $a \geq 1$ and $b \in \Z$
 
$$ \nu ^a : Gr^{{\cal W}'}_{b + a} Gr^{W^0}_{b} \Psi^u_f {\cal L}  \simeq Gr^{{\cal 
W}'}_{b - a} Gr^{W^0}_{b} \Psi^u_f {\cal L}$$ 
follows from the corresponding relation for $ Gr^{W^0}_b L $ as a $VHS$
and this concludes that: 
 ${\cal W}'$ is the relative monodromy filtration  with respect to $  W^0 $ on $ \Psi^u_f 
{\cal L}$.
This isomorphism is equivalent to

\smallskip

{\em Corollary: i) For all $a \geq 1, \; Gr^{{\cal W}'}_{b+a}Gr^{W^0}_{b} ker  \nu^a = 
s(Gr^{{\cal W}'}_{b+a+2p-1}Gr^{W^0}_{b}{\Omega}^* L [p], \eta)_{-a <p \leq 0} \cong 0$

ii) for $a > 0$, we have; $ Gr^{{\cal W}'}_{b + a} Gr^{W^0}_{b} \Psi^u_f {\cal L}  \simeq 
s(Gr^{{\cal W}'}_{b + a + 2p -1} Gr^{W^0}_{b}{\Omega}^* L [p], \eta))_{p  \leq -a} \simeq 
$

\noi $ \oplus_{a+ 2p -1 \leq -a-1}
Gr^{{\cal W}'}_{b + a + 2p -1} Gr^{W^0}_{b}{\Omega}^* L [p]$ where only a finite number 
of $p$ give non zero terms.

iii) for $t > v \in \Z, \; Gr^{{\cal W}'}_{r}(W^0_t/W^0_v)\Psi^u_f {\cal L}$ decomposes 
into a direct sum of intermediate extensions of 
 $VHS$ of weight $r$ graded polarised with respect to an induced filtration by $W^0$.

iv)for a proper morphism $f$, the complex 
$(\Psi^u_f {\cal L},{\cal W}', W^0, F)$ is a 
filtered mixed Hodge complex in the sense that for all $b < a $, 
$(W^0_a/W^0_b) \Psi^u_f {\cal L},{\cal W}', F)$ is 
a mixed Hodge complex.}

\smallskip

The proof is similar to the case of $VHS$ and uses the remark above for  successive  
extensions of $Gr^{W^0}_i L$. Notice the adjustment of the weight to $r$ in (iii)  since 
we take the sum for $r = a + b $ with $ a + 2p -1 \leq 0 $ that is the case of shift by $ 
 1 $ of the weight $ s + 1 $ of $(Gr^{{\cal W}'}_{s} {\Omega}^* Gr^{W^0}_{u}L$ for $ s 
\leq u $ which is compensated by $-1$ in the formula $r + 2p - 1$, while $2p$ is 
compensated by $p$ in the formula for $F$.
\bigskip

\bigskip

ACKNOWLEDGMENTS. I would like to thank J.P. Serre for accepting the note [14] and 
acnowledge the valuable discussions I  had at that time with J.L. Verdier when I started 
to work on the first preprints .

        D\'epart. de Math., Univ. de Nantes, 
        Laboratoire de  Math. Jean Leray
        CNRS - UMR 6629
        44072 Nantes Cedex 03, FRANCE. 

        e-mail: elzein@math.univ-nantes.fr

\end{document}